\documentclass[12pt,reqno]{amsart}
\usepackage[T1]{fontenc}
\usepackage[utf8]{inputenc}
\usepackage{amssymb, amsthm, amsxtra, epsfig, amsmath, parskip}
\usepackage{hyperref}
\usepackage{booktabs}
\usepackage{dcolumn}
\usepackage{color}
\usepackage{subfigure}

\newcommand{\ind}{{1\hspace{-1mm}{\rm I}}}
\newcommand{\N}{\mathbb{N}}

\newcommand{\R}{\mathbb{R}}
\newcommand{\Var}{\text{Var}}
\newcommand{\Cov}{\text{Cov}}
\newcommand{\E}{\mathbb{E}}

\newcommand{\Prob}{\mathbb{P}}
\newcommand{\Z}{\mathbb{Z}}

\newtheorem{theorem}{Theorem}[section]
\newtheorem{lemma}[theorem]{Lemma}

\newtheorem{corollary}[theorem]{Corollary}

\theoremstyle{definition}
\newtheorem{definition}[theorem]{Definition}
\newtheorem{example}[theorem]{Example}

\theoremstyle{remark}
\newtheorem{remark}[theorem]{Remark}

\begin{document}

\sloppy
\title[Simulation of infinitely divisible random fields]{Simulation of infinitely divisible random fields}

\author{Wolfgang Karcher}
\address{Wolfgang Karcher, Ulm University, Institute of Stochastics, Helmholtzstr. 18, 89081 Ulm, Germany}
\email{wolfgang.karcher\@@{}uni-ulm.de}
\author{Hans-Peter Scheffler}
\address{Hans-Peter Scheffler, University of Siegen, Fachbereich 6, Mathematik, Emmy-Noether-Campus, Walter-Flex-Str. 3, 57068 Siegen, Germany}
\email{scheffler\@@{}mathematik.uni-siegen.de}
\author{Evgeny Spodarev}
\address{Evgeny Spodarev, Ulm University, Institute of Stochastics, Helmholtzstr. 18, 89081 Ulm, Germany}
\email{evgeny.spodarev\@@{}uni-ulm.de}

\date{14 October 2009}

\begin{abstract}
Two methods to approximate infinitely divisible random fields are presented. The methods are based on approximating the kernel function in the spectral representation of such fields, leading to numerical integration of the respective integrals. Error bounds for the approximation error are derived and the approximations are used to simulate certain classes of infinitely divisible random fields.
\end{abstract}


\maketitle

\baselineskip=18pt

\section{Introduction}
\label{intro}

In many cases, the normal distribution is a reasonable model for real phenomena. If one considers the cumulative outcome of a great amount of influence factors, the normal distribution assumption can be justified by the Central Limit Theorem which states that the sum of a large number of independent and identically distributed random variables can be approximated by a normal distribution if the variance of these variables is finite. However, many real phenomena exhibit rather heavy tails. Stable distributions remedy this drawback by still being the limit distribution of a sum of independent and identically distributed random variables, but allowing for an infinite variance and heavy tails.

Stable distributions are a prominent example of the class of infinitely divisible distributions which we particularly concentrate on in this paper. Infinitely divisible distributions are distributions whose probability measure $\Prob$ is equal to the $n$-fold convolution of a probability measure $\Prob_n$ for any positive integer $n$. The class of infinitely divisible distributions comprises further well-known examples such as the Poisson, geometric, negative binomial, exponential, and gamma distribution, see~\cite{Sat05}, p.~(page) 21. These distributions are widely used in practice, for instance in finance to model the returns of stocks or in insurance to model the claim amounts and the number of claims of an insurance portfolio.

In order to include time dependencies or the spatial structure of real phenomena, random processes may be an appropriate model. If the dimension of the index set of the random process is greater than one, random processes are also called random fields. An example of infinitely divisible processes are L\'evy processes which have been extensively studied in the literature.

In this paper, we consider random fields that can be represented as a stochastic integral of a deterministic kernel as integrand and an infinitely divisible random measure as integrator. The kernel basically determines the dependence structure, whereas the infinitely divisible random measure inhibits the probabilistic characteristics of the random field. As already noted by \cite{CLL08}, practitioners have to try a variety of kernels and infinitely divisible random measures to find the model that best fits their needs.

Once the model is fixed, it is desirable to be able to perform simulations of the considered random field. There are several papers that are devoted to this problem. In \cite{BS08}, \cite{ST04} and \cite{WMZ04}, the fast Fourier transform is used for the simulation of linear fractional stable processes, whereas in \cite{D01}, a wavelet representation of a certain type of fractional stable processes was applied to simulate sample paths. Furthermore, \cite{CLL08} gives a general framework for the simulation of fractional fields.

In this paper, we consider infinitely divisible random fields for which the kernel functions are assumed to be Hölder-continuous or bounded which is a less restrictive assumption. Based on the respective assumption, we derive estimates for the approximation error when the kernel functions are approximated by step functions or by certain truncated wavelet series. The approximation allows for simulation since the integral representation of the random field reduces to a finite sum of random variables in this case.

In Section \ref{sec:approx}, we present the main results for the approximation error which is made when the kernel function is replaced by a step function or a truncated wavelet series. Section \ref{sec:simstud} is devoted to a brief simulation study where we apply the derived formulas for the approximation error to the simulation of two particular stable random fields. Finally, in Section \ref{sec:summary} we comment on the simulation results and the methods discussed in Section \ref{sec:approx}.

\section{Infinitely divisible random fields admitting an integral representation}
\label{sec:integral_represenation}

Let $\Lambda$ be an infinitely divisible random measure with control measure $\lambda$, cf.~(compare)~\cite{RR89}, pp.~(page and the following) 455 and \cite{JW94}, pp.~75. Let $f_t:\R^d \to \R$, $d \geq 1$, be $\Lambda$-integrable for all $t \in \R^q$, $q \geq 1$, that is there exists a sequence of simple functions $\{\tilde{f}_t^{(n)}\}_{n \in \N}$, $\tilde{f}_t^{(n)}:\R^d \to \R$ for all $t \in \R^q$, such that
\begin{enumerate}
 \item $\tilde{f}_t^{(n)} \to f_t \quad \lambda-\text{a.e.},$
 \item for every Borel set $B \subset \R^d$, the sequence $\{\int\limits_{B} \tilde{f}_t^{(n)}(x) \Lambda(dx)\}_{n \in \N}$ converges in probability.
\end{enumerate}
For each $t \in \R^q$, we define
$$ \int\limits_{\R^d} f_t(x) \Lambda(dx) := \underset{n \to \infty}{\text{plim}} \int\limits_{\R^d} \tilde{f}_t^{(n)}(x) \Lambda(dx),$$
cf.~\cite{RR89}, p.~460, where $\underset{n \to \infty}{\text{plim}}$ means convergence in probability, and consider random fields of the form
\begin{equation}
 X(t) = \int\limits_{\R^d} f_t(x) \Lambda(dx), \quad t \in \R^q. \label{eq:IDfield}
\end{equation}

\begin{remark}
 In \cite{RR89}, it is shown that $X(t)$ is infinitely divisible for all $t \in \R^q$, cf.~Theorem 2.7, p.~461 and the L\'evy form of the characteristic function of an inifinitely divisible random variable, p.~456. More generally, it can be shown that the random vector $(X(t_1),...,X(t_n))$ is infinitely divisible for all $t_1,...,t_n \in \R^q$ and $n \in \N$, see~\cite{KSS09c}. Therefore, any random field of the form (\ref{eq:IDfield}) is an infinitely divisible random field.
\end{remark}

\begin{example}\label{ex:rf_stable}
 Let $0<\alpha\leq 2$, $M$ be an (independently scattered) $\alpha$-stable random measure on $\R^d$ with control measure $m$ and skewness intensity $\beta$, see~\cite{ST94}, pp.~118. Furthermore, we assume that $f_t \in L^\alpha(\R^d,m)$ if $\alpha \neq 1$ and $f_t \in\{f \in L^1(\R^d,m): \int_{\R^d} \vert f(x) \ln \vert f(x)\vert \vert m(dx) < \infty\}$ if $\alpha = 1$ for all $t \in \R^q$. We denote the set of all functions $f_t$ satisfying these conditions by $\mathcal{F}$. Then
$$X(t) = \int\limits_{\R^d} f_t(x) M(dx), \quad t \in \R^q, $$
is an $\alpha$-stable random field.
\end{example}

\begin{example}\label{ex:rf_poisson}
 Let $\Phi$ be a Poisson random measure on $\R^d$ with intensity measure $\Theta$, see~\cite{SKM95}, p.~42. Furthermore, we assume that $f_t$ is a measurable function on $\R^d$ for each $t \in \R^q$. Then we can consider the shot noise field
 $$X(t) = \int\limits_{\R^d} f_t(x) \Phi(dx), \quad t \in \R^q, $$
 which can be written as
 \begin{equation}
  X(t) = \sum\limits_{x \in \Psi} f_t(x), \label{eq:poisson}
 \end{equation}
 where $\Psi$ is the support set of the Poisson random measure, cf.~\cite{SKM95}, p.~101.
\end{example}

\begin{example}
 Let $Q$ be a Poisson random measure on $(0,\infty) \times \R\setminus \{0\}$ with intensity measure $\mu \times \nu$. Here, $\mu$ is the Lebesgue measure and $\nu$ is the L\'evy measure, cf.~\cite{Pro05} and \cite{GMP07}. Let $G$ be a Gaussian ($2$-stable) random measure with Lebesgue control measure and skewness intensity $\beta \equiv 0$. Then
 $$X(t) = \int\limits_{\R^2} x \ind\{0 \leq s \leq t\} Q(ds,dx) - t \int\limits_{|x|<1} x \nu(dx) + \gamma t + \int\limits_{\R} \ind\{0\leq x \leq t\} G(dx)$$
 is the L\'evy process with L\'evy measure $\nu$, Gaussian part $G$ and drift $\gamma$, where
 $$ \gamma = \E\left( X_1 - \int_{|x| \geq 1} x \nu(dx) \right). $$
\end{example}

We now consider the cumulant function $C_{\Lambda(A)}(t) = \ln(\E e^{it\Lambda(A)})$ of $\Lambda(A)$ for a set $A$ in the $\delta$-ring $\mathcal{A}$ of bounded Borel subsets of $\R^d$ which is given by the L\'evy-Khintchine representation
$$C_{\Lambda(A)}(v) = iva(A) - \frac{1}{2}v^2b(A) + \int_\R (e^{ivr}-1-ivr\ind_{[-1,1]}(r))U(dr,A),$$
where $a$ is a $\sigma$-additive set function on $\mathcal{A}$, $b$ is a measure on the Borel $\sigma$-algebra $\mathcal{B}(\R^d)$, and $U(dr,A)$ is a measure on $\mathcal{B}(\R^d)$ for fixed $dr$ and a L\'evy measure on $\mathcal{B}(\R)$ for each fixed $A \in \mathcal{B}(\R^d)$, that is $U(\{0\},A) = 0$ and $\int_\R \min\{1,r^2\} U(dr,A) < \infty$, cf.~\cite{HPVJ08}, p.~605. The measure $U$ is referred to as the \textit{generalized L\'evy measure} and $(a,b,U)$ is called \textit{characteristic triplet}. The control measure $\lambda$ can be written as
\begin{equation}
 \lambda(A) = |a|(A) + b(A) + \int_\R \min\{1,r^2\} U(dr,A), \quad A \in \mathcal{A}, \label{eq:control_measure}
\end{equation}
where $|a|=a^+ + a^-$, see~\cite{RR89}, p.~456. Furthermore, $a$ and $b$ are absolutely continuous with respect to $\lambda$ and we have the formulas
\begin{eqnarray*}
 a(d\eta) = \tilde{a}(\eta) \lambda(d\eta), \quad
 b(d\eta) = \tilde{b}(\eta) \lambda(d\eta), \quad
 U(dr,d\eta) = V(dr,\eta) \lambda(d\eta),
\end{eqnarray*}
where $V(dr,\eta)$ is a L\'evy measure for fixed $\eta$, cf.~\cite{RR89}, p.~457.

We now introduce a so-called spot variable $L'(\eta)$ with cumulant function
$$C_{L'(\eta)}(v) = iv\tilde{a}(\eta) - \frac{1}{2} v^2 \tilde{b}(\eta) + \int_\R(e^{ivr}-1-ivr\ind_{[-1,1]}(r))V(dr,\eta)$$
with
\begin{eqnarray}
 \E(L'(\eta)) &=& \tilde{a}(\eta) + \int_{[-1,1]^C}r V(dr,\eta), \label{eq:evLprime} \\
 \Var(L'(\eta)) &=& \tilde{b}(\eta) + \int_\R r^2 V(dr,\eta) \label{eq:varLprime}
\end{eqnarray}
if $\E(L'(\eta))$ and $\Var(L'(\eta))$ exist. In \cite{HPVJ08}, p.~607, it is shown that
$$C_{X(t)}(v) = \int_{\R^d} C_{L'(\eta)}(vf_t(\eta))\lambda(d\eta).$$
We can use the cumulant function of $X(t)$ to obtain the second moment of $X(t)$ (in the case it exists):
\begin{equation}
 \E\left(X(t)^2\right) = \int_{\R^d} f_t^2(y) \Var(L'(y)) \lambda(dy) + \left( \int_{\R^d} f_t(y) \E(L'(y)) \lambda(dy) \right)^2. \label{eq:second_moment}
\end{equation}

As noted by \cite{HPVJ08}, for modelling purposes it is no resriction if we only consider characteristic triplets $(a,b,U)$ of the form
$$a(d\eta) = \tilde{a}_\nu(\eta)\nu(d\eta), \quad b(d\eta) = \tilde{b}_\nu(\eta)\nu(d\eta), \quad U(dr,d\eta) = V_\nu(dr,d\eta)\nu(d\eta),$$
where $\nu$ is a nonnegative measure on $\mathcal{B}(\R^d)$, $\tilde{a}_\nu:\R^d \to \R$ and $\tilde{b}_\nu:\R^d \to [0,\infty)$ are measurable functions, and $V_\nu(dr,\eta)$ is a L\'evy measure for fixed $\eta$.

\begin{example}\label{ex:rf_gamma}
 We choose the Lebesgue measure for $\nu$ and consider the characteristic triplet $(a,0,U)$ with
 \begin{eqnarray*}
  U(dr,d\eta) &=& V(dr,\eta) d\eta = \ind_{(0,\infty)}(r) \frac{1}{r} e^{-\theta r} dr d\eta \\
  a(d\eta) &=& \tilde{a}(\eta) d\eta = \frac{1}{\theta}\left(1-e^{-\theta}\right) d\eta
 \end{eqnarray*}
with $\theta \in (0,\infty)$. Then, by using (\ref{eq:evLprime}) and (\ref{eq:varLprime}), we get
\begin{eqnarray*}
 \E(L'(\eta)) = \frac{1}{\theta}, \quad 
 \Var(L'(\eta)) = \frac{1}{\theta^2},
\end{eqnarray*}
and by (\ref{eq:control_measure}), the control measure is proportional to the Lebesgue measure with
\begin{equation*}
 \lambda(d\eta) = \left(\frac{1+\theta-2\theta e^{-\theta}-e^{-\theta}}{\theta^2}+\int_1^\infty \frac{1}{r} e^{-\theta r} dr \right) d\eta.
\end{equation*}

\end{example}

\section{Approximation of infinitely divisible random fields}
\label{sec:approx}

We now restrict our setting to the observation window $[-T,T]^q$ with $T>0$ such that

$$X(t) = \int\limits_{\R^d} f_t(x) \Lambda(dx), \quad t \in [-T,T]^q. $$

We denote by $supp(f_t)$ the support of $f_t$ for each $t \in [-T,T]^q$ and assume that
$$ \bigcup\limits_{t \in [-T,T]^q} supp(f_t) \subset [-A,A]^d$$
for an $A > 0$. Then $X(\cdot)$ can be written as
$$X(t) = \int\limits_{[-A,A]^d} f_t(x) \Lambda(dx), \quad t \in [-T,T]^q. $$

Our goal is to approximate sample paths of $X$ for a variety of kernel functions $f_t$, $t \in [-T,T]^q$. The idea is to approximate the kernel functions $f_t$ appropriately such that the approximations $\tilde{f}_t^{(n)}$ are of the form
$$\tilde{f}_t^{(n)} = \sum_{i=1}^{m(n)} a_i g_{t,i}, \quad t \in [-T,T]^q,$$
where $m(n) \in \N$, $a_i \in \R$ and $g_{t,i}:\R^d \to \R$ is $\Lambda$-integrable.

Due to the linearity of the stochastic integral, we get as an approximation $\tilde{X}^{(n)}$ of $X$
$$ \tilde{X}^{(n)}(t) = \int_{[-A,A]^d} \tilde{f}_t^{(n)}(x) \Lambda(dx) = \sum_{i=1}^{m(n)} a_i \int_{[-A,A]^d} g_{t,i}(x) \Lambda(dx), \quad t \in [-T,T]^q.$$

If the $g_{t,i}$ are simple functions such that
$$\int_{[-A,A]^d} g_{t,i}(x) \Lambda(dx) = \sum_{j=1}^{l} g_{t,i}(x_j) \Lambda(\Delta_j), \quad i=1,...,m(n), \quad t \in [-T,T]^q,$$
for some $x_j \in [-A,A]^d$, $l \in \N$ and a partition $\{\Delta_j\}_{j=1}^{l}$ of $[-A,A]^d$, then
$$ \tilde{X}^{(n)}(t) = \sum_{i=1}^{m(n)}\sum_{j=1}^{l} a_i g_{t,i}(x_j) \Lambda(\Delta_j)$$
which can be simulated if $\Lambda(\Delta_j)$, $j=1,...,l$, can be simulated.

\begin{example}
 Let $\Lambda = M$ be an $\alpha$-stable random measure with Lebesgue control measure and constant skewness intensity $\beta$. Then
 $$M(\Delta_j) \sim S_\alpha(\vert \Delta_j \vert^{1/\alpha}, \beta, 0), \quad j=1,...,l,$$
 cf.~\cite{ST94}, p.~119, where $\vert \Delta_j \vert$ is the volume of $\Delta_j$ and $S_\alpha(\sigma,\beta,0)$ denotes the stable distribution with stable index $\alpha$, scale parameter $\sigma$, skewness parameter $\beta$ and location parameter $0$. Furthermore, $M(\Delta_j)$, $j=1,...,l$, are independent since $M$ is an independently scattered random measure. A method to simulate $\alpha$-stable random variables can be found in \cite{CMS76}.
\end{example}

\begin{example}
 Let $\Lambda = \Phi$ be a Poisson random measure with intensity measure $\Theta$. Then
 $$\Lambda(\Delta_j) \sim Poi(\Theta(\Delta_j)), $$
 where $Poi(\Theta(\Delta_j))$ denotes the Poisson distribution with mean $\Theta(\Delta_j)$. We note that simulating sample paths of $X$ by $\tilde{X}^{(n)}$ is not efficient since one can directly exploit the structure of $X$ and use
 $$ X(t) = \sum_{x \in \Psi} f_t(x),$$
 cf.~equation (\ref{eq:poisson}) in Example \ref{ex:rf_poisson}, p.~\pageref{ex:rf_poisson}.
\end{example}

\begin{example}
 Let $\Lambda_1=Q$ be a Poisson random measure on $(0,\infty)\times \R\setminus\{0\}$ with intensity measure $\mu \times \nu$, where $\mu$ is the Lebesgue measure and $\nu$ is the L\'evy measure. Let $\Lambda_2=G$ be a Gaussian random measure with Lebesgue control measure and skewness intensity $\beta \equiv 0$. We first approximate the L\'evy process
$$X(t) = \int\limits_{\R^2} x \ind\{0 \leq s \leq t\} Q(ds,dx) - t \int\limits_{|x|<1} x \nu(dx) + \gamma t + \int\limits_{\R} \ind\{0\leq x \leq t\} G(dx)$$
by
\begin{eqnarray*}
 X_K(t) &=& \int\limits_{-K}^K\int\limits_0^t x Q(ds,dx) - t \int\limits_{|x|<1} x \nu(dx) + \gamma t + \int\limits_{0}^t G(dx)\\
	&=& \int\limits_{-K}^K\int\limits_0^t x Q(ds,dx) - t \int\limits_{|x|<1} x \nu(dx) + \gamma t + G([0,t])
\end{eqnarray*}
for some $K>0$. We approximate $f(x) = x$ by using a linear combination of some simple functions $g_i$, $i=1,...,m(n)$, $m(n) \in \N$ for all $n \in \N$ and get for a partition $\{\Delta_j\}_{j=1}^l$, $l \in \N$, of $[-K,K]\times[0,t]$
$$\tilde{X}_K^{(n)}(t) = \sum_{i=1}^{m(n)}\sum_{j=1}^{l} a_i g_{i}(x_j)Q(\Delta_j) - t \int\limits_{|x|<1} x \nu(dx) + \gamma t + G([0,t]),$$
for some $a_1,...,a_{m(n)} \in \R$, where
\begin{eqnarray*}
 Q(\Delta_j) \sim Poi((\mu\times\nu)(\Delta_j)) \quad \text{and} \quad
 G([0,t]) \sim \mathcal{N}(\sqrt{t},0).
\end{eqnarray*}
\end{example}

\begin{example}
 We choose again the Lebesgue measure for $\nu$ and the characteristic triplet $(a,0,U)$ from Example \ref{ex:rf_gamma}, p.~\pageref{ex:rf_gamma}. Then
 $$\Lambda(\Delta_j) \sim \Gamma(|\Delta_j|,\theta),$$
 cf.~\cite{HPVJ08}, p.~608, where $|\Delta_j|$ is the Lebesgue measure of $\Delta_j$ and $\Gamma(|\Delta_j|,\theta)$ is the gamma distribution with probability density function
 $$f(x) = \frac{\theta^{|\Delta_j|}}{\Gamma(|\Delta_j|)} x^{|\Delta_j|-1} e^{-\theta x} \ind_{[0,\infty)}(x).$$
 In the last formula, $\Gamma(\cdot)$ denotes the gamma function. Again, $\Lambda(\Delta_j)$, $j=1,...,n$, are independent. Due to its distributional property, $\Lambda$ is called \textit{gamma L\'evy basis}.
\end{example}

\subsection{Measuring the approximation error}
\label{subsec:measure}
{\ }\newline
Approximating the random field $X$ with $\tilde{X}^{(n)}$ by taking an approximation $\tilde{f}_t^{(n)}$ of the kernel functions $f_t$ implicitely includes the assumption that $\tilde{X}^{(n)}$ is close to $X$ when $\tilde{f}_t^{(n)}$ is close to $f_t$. We use
$$ Err_s(X(t), \tilde{X}^{(n)}(t)) := \left\Vert f_t(x) - \tilde{f}_t^{(n)}(x) \right\Vert_{L^s} := \left(\int_{[-A,A]^d} |f_t(x) - \tilde{f}_t^{(n)}(x)|^s \lambda(dx) \right)^{1/s} $$
to measure the approximation quality of $\tilde{X}^{(n)}$ for an $s > 0$. In order to have existance of the above integral, we assume from now on that $f_t, \tilde{f}_t^{(n)} \in L^s([-A,A]^d,\lambda)$ for all $t \in [-T,T]^q$. The goal is then to find a set of functions $\{\tilde{f}_t^{(n)}\}_{t \in \R^d}$ such that $Err_s(X(t),\tilde{X}^{(n)}(t))$ is less than a predetermined critical value.

We see that the problem of approximating the random field $X$ reduces to an approximation problem of the corresponding kernel functions.

Let us now consider two special cases where $\Lambda$ is an $\alpha$-stable random measure and a Poisson random measure, respectively, to analyse the choice of the error measure.

\subsubsection{$\alpha$-stable random measures}
\label{approx_alpha_stable}
{\ }\newline
Assume that $0 < \alpha \leq 2$ and let $M$ be an $\alpha$-stable random measure with control measure $m$ and skewness intensity $\beta$. Furthermore, if $\alpha = 1$, assume additionally that $\beta(t) = 0$ for all $t \in [-T,T]^q$. Consider a set of functions $\{\tilde{f}_t^{(n)}\}_{t \in \R^d}$, where  $\tilde{f}_t^{(n)} \in \mathcal{F}$ (cf.~Example \ref{ex:rf_stable}, p.~\pageref{ex:rf_stable})  for all $t \in [-T,T]^q$ and $n \in \N$. The corresponding $\alpha$-stable random field is denoted by
$$ \tilde{X}^{(n)}(t) := \int\limits_{[-A,A]^d} \tilde{f}_t^{(n)}(x)M(dx), \quad t \in [-T,T]^q.$$

We know that $\tilde{X}^{(n)}(t)$ converges to $X(t)$ in probability if and only if \linebreak $\int_{\R^d} |f_t(x) - \tilde{f}_t^{(n)}(x)|^\alpha m(dx)$ converges to $0$ as $n$ goes to infinity, see~\cite{ST94}, p.~126. Therefore, we can use $\tilde{X}^{(n)}(t)$ as an approximation for $X(t)$ if $\tilde{f}_t^{(n)}$ approximates $f_t$ sufficiently well and $\tilde{X}^{(n)}(t)$ converges to $X(t)$ in probability if and only if
$$ Err_\alpha(X(t), \tilde{X}^{(n)}(t)) \to 0, \quad n \to \infty. $$
The choice of $ Err_\alpha(X(t), \tilde{X}^{(n)}(t))$ can be further justified as follows.

Since $X(t)$ and $\tilde{X}^{(n)}(t)$ are jointly $\alpha$-stable random variables for all $t \in [-T,T]^q$, the difference $X(t) - \tilde{X}^{(n)}(t)$ is also an $\alpha$-stable random variable. The scale parameter of $X(t) - \tilde{X}^{(n)}(t)$ is given by
$$\sigma_{X(t)-\tilde{X}^{(n)}(t)} = \left(\int_{[-A,A]^d} |f_t(x) - \tilde{f}_t^{(n)}(x)|^\alpha m(dx)\right)^{1/\alpha},$$
cf.~\cite{ST94}, p.~122, so that
$$ Err_\alpha(X(t), \tilde{X}^{(n)}(t)) = \sigma_{X(t)-\tilde{X}^{(n)}(t)}.$$

Furthermore, let us consider the quantity
$$\E |X(t)-\tilde{X}^{(n)}(t)|^p, \quad 0 < p < \alpha,$$
that is the mean error between $X(t)$ and $\tilde{X}^{(n)}(t)$ in the $L^p$-sense.

Since $X(t)-\tilde{X}^{(n)}(t)$ is an $\alpha$-stable random variable, we have
$$\E |X(t)-\tilde{X}^{(n)}(t)|^p < \infty, \quad 0 < p < \alpha$$
and
$$\E |X(t)-\tilde{X}^{(n)}(t)|^p = \infty, \quad p \geq \alpha.$$
For $0 < p < \alpha$, $0<\alpha<2$ and $\alpha \neq 1$, this quantity can be written as
\begin{equation}
 \left( \E |X(t)-\tilde{X}^{(n)}(t)|^p \right)^{1/p} = c_{\alpha, \beta_t}(p) \cdot \sigma_{X(t)-\tilde{X}^{(n)}(t)}, \label{eq:ev}
\end{equation}
where
\begin{equation*}
 \left(c_{\alpha, \beta_t}(p)\right)^p = \frac{2^{p-1}\Gamma(1-\frac{p}{\alpha})}{p \int_0^\infty u^{-p-1} \sin^2 u\ du} \left(1+\beta_t^2 \tan^2 \frac{\alpha \pi}{2} \right)^{p/2\alpha} \cos \left(\frac{p}{\alpha} \arctan\left(\beta_t \tan \frac{\alpha \pi}{2}\right)\right)
\end{equation*}
and
$$\beta_t = \frac{\int_{[-A,A]^d} |f_t(x) - \tilde{f}_t^{(n)}(x)|^\alpha \text{sign}(f_t(x) - \tilde{f}_t^{(n)}(x)) \beta(x) dx}{\int_{[-A,A]^d} |f_t(x) - \tilde{f}_t^{(n)}(x)|^\alpha dx}. $$

In the case $\alpha = 1$, equation (\ref{eq:ev}) holds if $\beta_t = 0$, see~\cite{ST94}, p.~18.

We remind that $\beta(\cdot)$ is the skewness intensity of the $\alpha$-stable random measure $M$. The above implies that
$$Err_\alpha(X(t), \tilde{X}^{(n)}(t)) = \frac{1}{c_{\alpha, \beta_t}(p)} \cdot \left(\E |X(t) - \tilde{X}^{(n)}(t)|^p \right)^{1/p}, \quad 0 < p < \alpha.$$

Now assume that $\alpha = 1$ and $\beta(t) \neq 0$ for at least one $t \in [-T,T]^q$. In this case, we need to impose an additional condition on the kernel functions in order to guarantee the convergence of $\tilde{X}^{(n)}(t)$ to $X(t)$ in probability. Namely, $\tilde{X}^{(n)}(t)$ converges to $X(t)$ in probability if and only if
\begin{eqnarray*}
 &\int_{[-A,A]^d}\vert f_t(x) - \tilde{f}_t^{(n)}(x)\vert m(dx) \to 0, \quad n \to \infty
\end{eqnarray*}
and
\begin{eqnarray*}
 &\int_{[-A,A]^d} (f_t(x) - \tilde{f}_t^{(n)}(x)) \ln \vert f_t(x) - \tilde{f}_t^{(n)}(x)\vert \beta(x) m(dx) \to 0, \quad n \to \infty,&
\end{eqnarray*}
cf.~\cite{ST94}, p.~126. We use
$$ Err_{3/2}(X(t), \tilde{X}^{(n)}(t)) := \left(\int_{[-A,A]^d} |f_t(x) - \tilde{f}_t^{(n)}(x)|^{3/2} m(dx) \right)^{2/3}$$
to measure the approximation error.

By using the fact that $ \vert x \ln x \vert \leq \max\{\sqrt{x}, x \sqrt{x}\}$ for $x > 0$ and that $-1 \leq \beta(t) \leq 1$ for all $t \in [-T,T]^q$, we have
\begin{eqnarray*}
 &&\left\vert\int_{[-A,A]^d} (f_t(x) - \tilde{f}_t^{(n)}(x)) \ln \vert f_t(x) - \tilde{f}_t^{(n)}(x)\vert \beta(x) m(dx) \right\vert\\
 &\leq& \int_{[-A,A]^d} \left\vert f_t(x) - \tilde{f}_t^{(n)}(x)) \ln \vert f_t(x) - \tilde{f}_t^{(n)}(x)\vert \beta(x) \right\vert m(dx) \\
 &\leq& \int_{[-A,A]^d} \max\left\{\vert f_t(x) - \tilde{f}_t^{(n)}(x)\vert^{1/2}, \vert f_t(x) - \tilde{f}_t^{(n)}(x)\vert^{3/2} \right\} m(dx) \\
 &\leq& \int_{[-A,A]^d} \vert f_t(x) - \tilde{f}_t^{(n)}(x)\vert^{1/2} m(dx) + \int_{[-A,A]^d} \vert f_t(x) - \tilde{f}_t^{(n)}(x)\vert^{3/2} m(dx).
\end{eqnarray*}

Now, if $Err_{3/2}(X(t), \tilde{X}^{(n)}(t))$ tends to $0$ as $n$ goes to infinity, then by Lyapunov's inequality (see~\cite{Shi95}),
\begin{eqnarray*}
 &&\int_{[-A,A]^d} \vert f_t(x) - \tilde{f}_t^{(n)}(x)\vert^{3/2} m(dx) \to 0, \quad n \to \infty, \\
 &&\int_{[-A,A]^d} \vert f_t(x) - \tilde{f}_t^{(n)}(x)\vert^{1/2} m(dx) \to 0, \quad n \to \infty,
\end{eqnarray*}

such that
$$\int_{[-A,A]^d} (f_t(x) - \tilde{f}_t^{(n)}(x)) \ln \vert f_t(x) - \tilde{f}_t^{(n)}(x)\vert \beta(x) m(dx) \to 0, \quad n \to \infty.$$

Furthermore, once again by using Lyapunov's inequality,
$$\int_{[-A,A]^d} \vert f_t(x) - \tilde{f}_t^{(n)}(x)\vert m(dx) \to 0, \quad n \to \infty,$$
so that $\tilde{X}^{(n)}(t)$ converges to $X(t)$ in probability.

\begin{remark}
 If $Err(X(t), \tilde{X}^{(n)}(t))$ tends to $0$ as $n$ goes to infinity for all $t \in [-T,T]^q$, then
 $$ \{\tilde{X}^{(n)}(t)\}_{t \in [-T,T]^q} \stackrel{f.d.}{\rightarrow} \{X(t)\}_{t \in [-T,T]^q}, $$
 where $\stackrel{f.d.}{\to}$ denotes convergence in distributions of all finite dimensional marginals.
\end{remark}

\begin{proof}
 We first state a lemma which is an implication of the inequality 
$$(x+y)^p \leq x^p+y^p, \quad 0<p\leq 1, \quad x,y>0.$$

\begin{lemma}\label{lemma:inequality}
 Let $0<p \leq 1$ and $f_i \in L^p(\R^d)$ for $i=1,...,n$ with $n \in \N$. Then
 $$ \left\Vert \sum_{i=1}^n f_i \right\Vert_{L^p}^p \leq \sum_{i=1}^n \Vert f_i \Vert_{L^p}^p.$$
\end{lemma}
 Now fix any $t_1,...,t_m \in [-T,T]^q$ and $\lambda_1,...,\lambda_m \in \R$. If $0 < \alpha < 1$, we get
 \begin{eqnarray*}
  &&\left(Err_\alpha\left(\sum\limits_{j=1}^m \lambda_j X(t_j), \sum\limits_{j=1}^m \lambda_j \tilde{X}^{(n)}(t_j)\right)\right)^\alpha \\
   &=& \left\Vert \sum\limits_{j=1}^m \lambda_j f_{t_j} - \sum\limits_{j=1}^m \lambda_j \tilde{f}^{(n)}_{t_j} \right\Vert_{L^\alpha}^\alpha
  \leq \sum\limits_{j=1}^m \vert \lambda_j \vert \left\Vert f_{t_j} - \tilde{f}^{(n)}_{t_j} \right\Vert_{L^\alpha}^\alpha \\
  &=& \sum\limits_{j=1}^m \vert \lambda_j \vert Err_\alpha(X(t_j), \tilde{X}^{(n)}(t_j))^\alpha \to 0, \quad n \to \infty.
 \end{eqnarray*}
 
 For $1 < \alpha \leq 2$, we can use Minkowski's inequality and get
 \begin{eqnarray*}
  &&Err_\alpha\left(\sum\limits_{j=1}^m \lambda_j X(t_j), \sum\limits_{j=1}^m \lambda_j \tilde{X}^{(n)}(t_j)\right) \\
   &=& \left\Vert \sum\limits_{j=1}^m \lambda_j f_{t_j} - \sum\limits_{j=1}^m \lambda_j \tilde{f}^{(n)}_{t_j} \right\Vert_{L^\alpha}
  \leq \sum\limits_{j=1}^m \vert \lambda_j \vert \left\Vert f_{t_j} - \tilde{f}^{(n)}_{t_j} \right\Vert_{L^\alpha} \\
  &=& \sum\limits_{j=1}^m \vert \lambda_j \vert Err_\alpha(X(t_j), \tilde{X}^{(n)}(t_j)) \to 0, \quad n \to \infty.
 \end{eqnarray*}
 
 Analogously for $\alpha = 1$,
$$Err_{3/2}\left(\sum\limits_{j=1}^m \lambda_j X(t_j), \sum\limits_{j=1}^m \lambda_j \tilde{X}^{(n)}(t_j)\right) \to 0,\quad n \to \infty.$$

Therefore, for any $\lambda_1,...,\lambda_m \in \R$ we have
$$\sum\limits_{j=1}^m \lambda_j X(t_j) - \sum\limits_{j=1}^m \lambda_j \tilde{X}^{(n)}(t_j) \to 0 $$
in probability which implies the convergence of all finite-dimensional marginal distributions.
\end{proof}

\subsubsection{Poisson random measures}
\label{approx_Poisson}
{\ }\newline
Let $\Phi$ be a Poisson random measure with intensity measure $\Theta$ and $\Psi$ be the random sequence of the Poisson point process corresponding to the Poisson random measure. Furthermore, assume that $f_t$ is measurable on $[-A,A]^d$ for each $t \in [-T,T]^q$. Then, by the Campbell theorem (see~\cite{SKM95}, p.~103), we have
\begin{eqnarray*}
 Err_1(X(t), \tilde{X}^{(n)}(t)) &=& \int\limits_{[-A,A]^d} \left\vert f_t(x) - \tilde{f}_t^{(n)}(x) \right\vert \Theta(dx) \\
 &=& \E \left(\ \int\limits_{[-A,A]^d} \left\vert f_t(x) - \tilde{f}_t^{(n)}(x) \right\vert \Phi(dx) \right) \\
 &=& \E \left( \sum_{x \in \Psi} \left\vert f_t(x) - \tilde{f}_t^{(n)}(x) \right\vert \right)
 \geq \E \left\vert \sum_{x \in \Psi} \left( f_t(x) - \tilde{f}_t^{(n)}(x) \right) \right\vert \\
 &=& \E \left\vert \int_{[-A,A]^d} f_t(x) \Phi(dx) - \int_{[-A,A]^d} \tilde{f}_t^{(n)}(x) \Phi(dx) \right\vert \\
 &=& \E \left\vert X(t) - \tilde{X}^{(n)}(t) \right\vert,
\end{eqnarray*}
that is we can control the mean error between $X(t)$ and $\tilde{X}^{(n)}(t)$ in the $L^1$-sense by finding error bounds for $Err_1(X(t), \tilde{X}^{(n)}(t))$.

\subsubsection{Exploiting the spot variable representation of the second moment of $X(t)$}
{\ }\newline
We assume that the second moment of the random field exists and recall formula (\ref{eq:second_moment}), p.~\pageref{eq:second_moment},
\begin{equation*}
 \E\left(X(t)^2\right) = \int_{[-A,A]^d} f_t^2(y) \Var(L'(y)) \lambda(dy) + \left( \int_{[-A,A]^d} f_t(y) \E(L'(y)) \lambda(dy) \right)^2
\end{equation*}
which implies
\begin{eqnarray*}
 \E(X(t)-\tilde{X}^{(n)}(t))^2 &=& \int\limits_{[-A,A]^d} (f_t(y)-\tilde{f}_t^{(n)}(y))^2 \Var(L'(y)) \lambda(dy) \\
 && + \left(\  \int\limits_{[-A,A]^d} (f_t(y)-\tilde{f}_t^{(n)}(y)) \E(L'(y)) \lambda(dy) \right)^2 \\
 &\leq& \int\limits_{[-A,A]^d} (f_t(y)-\tilde{f}_t^{(n)}(y))^2 \Var(L'(y)) \lambda(dy) \\
 && + \int\limits_{[-A,A]^d} (f_t(y) - \tilde{f}_t^{(n)}(y))^2 \lambda(dy) \cdot \int\limits_{[-A,A]^d} \left(\E(L'(y))\right)^2 \lambda(dy) \\
\end{eqnarray*}
where we used the Cauchy-Schwarz inequality in the last inequality. If
\begin{eqnarray*}
 \Var(L'(y)) \leq c_1 < \infty, \quad \forall y \in \R^d \quad \text{and} \quad 
 \int_{[-A,A]^d} \left(\E(L'(y))\right)^2 \lambda(dy) := c_2 < \infty,
\end{eqnarray*}
then we get
\begin{equation*}
 \left(\E(X(t) - \tilde{X}^{(n)}(t))^2\right)^{1/2} \leq (c_1 + c_2)^{1/2} \Vert f_t - g_t \Vert_{L^2} = (c_1 + c_2)^{1/2}Err_2(X(t), \tilde{X}^{(n)}(t)).
\end{equation*}

\begin{example}
 We choose again the Lebesgue measure for $\nu$ and the characteristic triplet $(a,0,U)$ from Example \ref{ex:rf_gamma}, p.~\pageref{ex:rf_gamma}. We have
 \begin{eqnarray*}
  \Var(L'(y)) &=& \frac{1}{\theta^2} =: c_1 \quad \text{and} \quad
  \int\limits_{[-A,A]^d} \left(\E(L'(y))\right)^2 \lambda(dy) = \frac{(2A)^d}{\theta^2} =: c_2,
 \end{eqnarray*}
 such that
 $$\left(\E(X(t) - \tilde{X}^{(n)}(t))^2\right)^{1/2} \leq \frac{1}{\theta}\left(1 + (2A)^d\right)^{1/2} Err_2(X(t), \tilde{X}^{(n)}(t)).$$
\end{example}

\subsection{Step function approximation}
\label{subsec:step}
{\ }\newline
For any natural number $n \geq 1$ and $k=(k_1,...,k_d) \in \Z^d$ with $-n \leq k_1,...,k_d < n$, let
\begin{eqnarray*}
 \xi_{k} &=& \left(k_1 \frac{A}{n},\cdots, k_d \frac{A}{n}\right), \\
 \Delta_{k} &=& \left[ k_1 \frac{A}{n}, (k_1+1) \frac{A}{n}\right) \times \cdots \times \left[ k_d \frac{A}{n}, (k_d+1) \frac{A}{n}\right).
\end{eqnarray*}

We define the step function
\begin{equation*}
 \tilde{f}_t^{(n)}(x) := \sum\limits_{|k| \leq n} f_t(\xi_k) \ind_{\Delta_{k}}(x)
\end{equation*}
to approximate $f_t$, where $|k| \leq n$ is meant to be componentwise, i.~e. $-n\leq k_i<n$ for $i=1,...,d$. Then we have
\begin{equation}
\tilde{X}^{(n)}(t) = \int\limits_{[-A,A]^d} \tilde{f}_t^{(n)}(x) \Lambda(dx) = \sum\limits_{|k| \leq n} f_t(\xi_k) \Lambda(\Delta_{k}). \label{eq:rf_step_function}
\end{equation}

The following theorem provides error bounds for $Err_s(X(t),\tilde{X}^{(n)}(t))$ for H\"older-continuous functions $f_t$.

\begin{theorem} \label{th:hoelder}
 Assume that $0<s\leq 2$, the control measure $\lambda$ is the Lebesgue measure and the functions $f_t$ are Hölder-continuous for all $t \in [-T,T]^q$, i.~e.
 $$ |f_t(x) - f_t(y)| \leq C_t \cdot ||x - y||_2^{\gamma_t}, \quad x,y \in [-A,A]^d, \quad t \in [-T,T]^q$$
for some $0 < {\gamma_t} \leq 1$ and $C_t > 0$, where $\Vert \cdot \Vert_2$ denotes the Euclidean norm. Then for any $t \in [-T,T]^q$ we have for all $n \geq 1$ that
\begin{equation}
 Err_s(X(t),\tilde{X}^{(n)}(t)) \leq \left(\frac{2^d C_t d}{1+{\gamma_t} s}\right)^{1/s} A^{{\gamma_t}+d/s} \left(\frac{1}{n}\right)^{{\gamma_t}}.\label{eq:hoelder}
\end{equation}
\end{theorem}

\begin{proof}
Since
$$X(t) - \tilde{X}^{(n)}(t) = \int_{[-A,A]^d} \left( f_t(x) - \sum\limits_{|k| \leq n} f_t(\xi_k) \ind_{\Delta_{k}}(x) \right) \Lambda(dx),$$
we have
\begin{eqnarray*}
 \left(Err_s(X(t),\tilde{X}^{(n)}(t))\right)^s &=& \int\limits_{[-A,A]^d} \left|f_t(x) - \sum_{|k| \leq n} f_t(\xi_k) \ind_{\Delta_k}(x)\right|^s  dx \\
	&=& \int\limits_{[-A,A]^d} \left| \sum_{|k| \leq n} \left(f_t(x) - f_t(\xi_k)\right) \ind_{\Delta_k}(x) \right|^s  dx.
\end{eqnarray*}

For each $x \in [-A,A]^d$, there exists exactly one $\tilde{k} = \tilde{k}(x)$ with $\vert \tilde{k} \vert \leq n$ and $x \in \Delta_{\tilde{k}}$. Hence
\begin{eqnarray*}
 \left| \sum_{|k| \leq n} \left(f_t(x) - f_t(\xi_k)\right) \ind_{\Delta_k}(x) \right|^s 
 &=& \left|\left(f_t(x) - f_t(\xi_k)\right)\right|^s \ind_{\Delta_{\tilde{k}}}(x) \\
 &=& \sum_{|k| \leq n} \left|f_t(x) - f_t(\xi_k)\right|^s \ind_{\Delta_k}(x),
\end{eqnarray*}
which implies
\begin{eqnarray*}
 \left(Err_s(X(t),\tilde{X}^{(n)}(t))\right)^s &=& \int\limits_{[-A,A]^d} \left| \sum_{|k| \leq n} \left(f_t(x) - f_t(\xi_k)\right) \ind_{\Delta_k}(x) \right|^s  dx \\
	&=& \sum_{|k| \leq n} \hspace*{0.2cm} \int\limits_{\Delta_k} \left|\left(f_t(x) - f_t(\xi_k)\right)\right|^s dx \\
	&\leq& C_t \sum_{|k| \leq n} \int\limits_{\Delta_{k}} ||x-\xi_k||_2^{{\gamma_t} s} dx \\
	&=& C_t \sum_{|k| \leq n} \int\limits_0^{A/n}\cdots \int\limits_0^{A/n} \left(y_1^2 + \cdots + y_d^2\right)^{({\gamma_t} s)/2} dy_d \cdots dy_1.
\end{eqnarray*}

As $({\gamma_t} s)/2 \leq 1$, we have $\left(y_1^2+\cdots +y_d^2 \right)^{({\gamma_t} s)/2} \leq y_1^{{\gamma_t} s}+\cdots +y_d^{{\gamma_t} s}$ and hence
\begin{eqnarray}
&&\int\limits_0^{A/n}\cdots \int\limits_0^{A/n} \left(y_1^2 + \cdots + y_d^2\right)^{({\gamma_t} s)/2} dy_d \cdots dy_1 \nonumber\\
&\leq& d \left(\frac{A}{n}\right)^{d-1} \int\limits_0^{A/n} y_1^{{\gamma_t} s} dy_1 
 = \frac{d}{{\gamma_t} s+1} \left(\frac{A}{n}\right)^{d+{\gamma_t} s}. \label{eq:control}
\end{eqnarray}

Therefore, we get
$$ Err_s(X(t),\tilde{X}^{(n)}(t)) \leq \left(\frac{2^d C_t d}{1+{\gamma_t} s}\right)^{1/s} A^{{\gamma_t}+d/s} \left(\frac{1}{n}\right)^{{\gamma_t}}.$$ \qed
\end{proof}

\begin{remark}
 It suffices to consider a control measure $\lambda$ proportional to the Lebesgue measure (cf.~Example \ref{ex:rf_gamma}, p.~\ref{ex:rf_gamma}), that is
 $$\lambda(d\eta) = c \cdot d\eta, \quad c > 0.$$
 In this case, one has to multiply the upper bound in (\ref{eq:hoelder}) by $c^{1/2}$. If the control measure $\lambda$ is not the Lebesgue measure, then, in general, the integral
$$ \int\limits_0^{A/n}\cdots \int\limits_0^{A/n} \left(y_1^2 + \cdots + y_d^2\right)^{({\gamma_t} s)/2} \lambda(d(y_1,...,y_d))$$
in (\ref{eq:control}) cannot be calculated explicitely such that one would have to include it in the upper bound of the approximation error.
\end{remark}

\begin{remark}
 In the proof, one can estimate
 $$\int\limits_0^{A/n}\cdots \int\limits_0^{A/n} \left(y_1^2 + \cdots + y_d^2\right)^{({\gamma_t} s)/2} dy_d\cdots dy_1$$
 alternatively by
 \begin{eqnarray*}
  \int\limits_0^{A/n}\cdots \int\limits_0^{A/n} \Vert y \Vert_2^{\gamma_t s} dy_d\cdots dy_1
  \leq \frac{1}{4} \int\limits_{\Vert y \Vert_2 \leq \frac{A}{n}\sqrt{d}} \Vert y \Vert_2^{\gamma_t s} dy
 \end{eqnarray*}
 and calculate the last integral by using polar coordinates. This yields
 \begin{eqnarray*}
 Err_s(X(t),\tilde{X}^{(n)}(t))
 \leq \left(\frac{2^d C_t\frac{\pi}{2} d^{(\gamma_t s + d)/2}}{\gamma_t s +d}\right)^{1/s} A^{\gamma_t + d/s} \left(\frac{1}{n}\right)^{\gamma_t} D(d,s),
\end{eqnarray*}
where
$$ D(d,s) := \begin{cases}
  1, & d=2, \\
  2^{1/s}, & d=3, \\
  \pi^{1/s}, & d=4, \\
  \left(\pi^{d-3} \cdot \frac{\Gamma(3/2)}{\Gamma(d/2)}\right)^{1/s}, & d \geq 5 \text{ odd}, \\
  \left(\pi^{d-7/2} \cdot \frac{\Gamma(3/2)}{\Gamma((d-1)/2)}\right)^{1/s}, & d \geq 5 \text{ even}.
 \end{cases}$$
It is straightforward to show that for $d=2$, this estimate is worse than the one in Theorem \ref{th:hoelder} if
$$\frac{\pi}{2} d^{(\gamma_t s + d)/2 -1} \cdot D^s(d,s) \cdot \frac{1+\gamma_t s}{\gamma_t s+d} \geq 1.$$
This is illustrated in Figure \ref{fig:combinations}. For $d \geq 3$, however, the estimate in Theorem \ref{th:hoelder} always performs better.

\begin{figure}[ht]
\includegraphics[width=7cm]{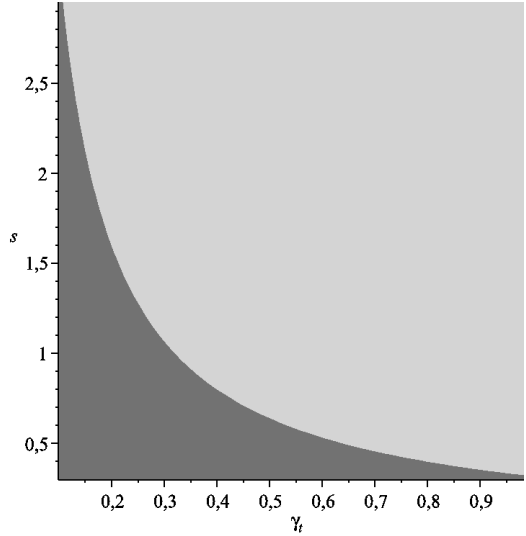}
\caption{Combinations of $(\gamma_t,s)$ for $d=2$. For all combinations in the dark grey area, the constant in Theorem \ref{th:hoelder} performs better, for all other combinations (in the light grey area) it is vice versa.}
\label{fig:combinations}
\end{figure}
\end{remark}

\begin{remark}
 Suppose that the conditions of Theorem \ref{th:hoelder} hold true. If the support of $f_t$ is not compact, we first need to estimate
 $$ X(t) = \int\limits_{\R^d} f_t(x) \Lambda(dx) $$
 by
 $$ X_K(t) = \int\limits_{[-K,K]^d} f_t(x) \Lambda(dx) $$
 For $K > 0$ large enough, the approximation error is small since
 $$ Err_s(X(t), X_K(t)) = \left(\ \int\limits_{\R^d\setminus [-K,K]^d} |f_t(x)|^s dx\right)^{1/s} \to 0, \quad K \to \infty $$
 and $f_t, \tilde{f}_t^{(n)} \in L^s([-A,A]^d,\lambda)$. Let $\varepsilon > 0$. If $1 \leq s \leq 2$, choose $K > 0$ such that $Err_s(X(t), X_K(t)) \leq \varepsilon /2$. We can apply Theorem \ref{th:hoelder} to $X_K(\cdot)$ such that 
$$Err_s(X_K(t),\tilde{X}_K^{(n)}(t)) \leq \frac{\varepsilon}{2}$$
for $n \in \N$ large enough. Then
\begin{eqnarray*}
 Err_s(X(t), \tilde{X}_K^{(n)}(t)) \leq Err_s(X(t),X_K(t)) + Err_s(X_K(t), \tilde{X}_K^{(n)}(t)) \leq \frac{\varepsilon}{2} + \frac{\varepsilon}{2} =\varepsilon.
\end{eqnarray*}
If $0 < s < 1$, choose $K > 0$ such that $Err_s(X(t), X_K(t)) \leq \varepsilon^s /2$. Again, we can apply Theorem \ref{th:hoelder} to $X_K(\cdot)$ such that 
$$Err_s(X_K(t),\tilde{X}_K^{(n)}(t)) \leq \frac{\varepsilon^s}{2}$$
for $n \in \N$ large enough. Then
\begin{eqnarray*}
 Err_s(X(t), \tilde{X}_K^{(n)}(t)) &\leq& \left(Err_s(X(t),X_K(t)) + Err_s(X_K(t), \tilde{X}_K^{(n)}(t))\right)^{1/s} \\
 	&\leq& \left(\frac{\varepsilon^s}{2} + \frac{\varepsilon^s}{2}\right)^{1/s} =\varepsilon.
\end{eqnarray*}
\end{remark}

\begin{remark}
 Theorem \ref{th:hoelder} provides a pointwise estimate of the approximation error for each $t \in [-T,T]^q$. We can obtain a uniform error bound as follows.

Assume that $\gamma := \inf\limits_{t \in [-T,T]^q} \gamma_t > 0$. Then for each $t \in [-T,T]^q$, $f_t$ is H\"older-continuous with parameters $\gamma$ and some constant $C_t^* > 0$. Set $C := \sup\limits_{t \in [,T,T]^d} C_t^*$. Then $Err_s(X(t),\tilde{X}^{(n)}(t))$ can be estimated by (\ref{eq:hoelder}) with $C_t$ and $\gamma_t$ replaced by $C$ and $\gamma$.

One can also consider the integrated error
$$Err_s(X,\tilde{X}^{(n)}) := \int\limits_{[-T,T]^q} Err_s(X(t),\tilde{X}^{(n)}(t)) dt$$
and multiply the error bound by $(2T)^q$.
\end{remark}

\begin{remark}
 Assume that $0 < s \leq 2$ and the functions $f_t$ are differentiable with $||\nabla f_t(x)||_2 \leq C_t$ for all $x \in [-A,A]^d$ and $t \in [-T,T]^q$. Then for any $t \in [-T,T]^q$, $(\ref{eq:hoelder})$ holds for all $n \geq 1$ with ${\gamma_t} = 1$, that is
\begin{equation*}
 Err_s(X(t),\tilde{X}^{(n)}(t)) \leq \left(\frac{2^d C_t d}{1+s}\right)^{1/s} A^{1+d/s} \cdot \frac{1}{n}
\end{equation*}
since $f_t$ is Hölder-continuous with $C_t$ and ${\gamma_t}=1$.
\end{remark}

\subsection{Approximation by wavelet series}
\label{subsec:wavelets}

\subsubsection{Series representation of kernel functions}
\label{subsubsec:series}
\vspace*{-0.3cm}
{\ }\newline
Let $s > 0$, $f_t \in L^s(\R^d)$, $t \in \R^q$, and let $\{\xi_i\}_{i \in I}$ be a basis for $L^s(\R^d)$, where $I$ is an index set. Then $f_t$ can be represented as

\begin{equation}
 f_t = \sum\limits_{i \in I} a_i \cdot \xi_i \label{eq:series}
\end{equation}

for certain constants $a_i \in \R$. In order to approximate $f_t$, one can truncate (\ref{eq:series}) such that it consists only of a finite number of summands. In \cite{BS08}, the trigonometric system is used to approximate the kernel function of certain stable random fields. In this paper, we will go another way and analyse whether a wavelet system may be appropriate for the simulation of random fields with an infinitely divisible random measure as integrator.

\subsubsection{Haar wavelets}
\label{subsubsec:introduction}
{\ }\newline
In this paper, we use the so-called \textit{Haar basis} to approximate the kernel functions. For a detailed introduction into wavelets, see for example \cite{Urb08}, \cite{DL92} and \cite{DeV98}.

\begin{definition}
 Consider the function
 $$ \varphi^{\text{Haar}}(x) := \frac{1}{(2A)^{1/2}} \cdot \ind_{[-A,A]}(x), \quad x \in \R,$$
 and the corresponding \textit{mother wavelet} defined by
 $$ \Psi^{\text{Haar}}(x) := \varphi^{\text{Haar}}(2x+A) - \varphi^{\text{Haar}}(2x-A), \quad x \in \R.$$
 The resulting basis $\mathcal{H}:=\{\varphi^{\text{Haar}}\} \cup \{\Psi_{j-2^k,k}^{\text{Haar}}\}_{\substack{k \in \N_0, \\ 2^k \leq j \leq 2^{k+1}-1}}$ is called Haar basis for $L^{2}([-A,A])$, where $\Psi_{j,k}^{\text{Haar}}(x) := 2^{k/2} \Psi^{\text{Haar}} (2^k (x+A) - (1+2\cdot j)A)$.
\end{definition}

Let $\Psi^0 := \varphi^{\text{Haar}}$, $\Psi^1 := \Psi^{\text{Haar}}$ and $E$ be the set of nonzero vertices of the unit cube $[0,1]^d$. Consider the multivariate functions $\Psi^e$, $e=(e_1,...,e_d) \in E$, defined by
$$\Psi^e(x_1,...,x_d) := \Psi^{e_1}(x_1)\cdots \Psi^{e_d}(x_d), \quad x \in \R^d.$$
Let $x=(x_1,...,x_d)$ and $a=(A,...,A)^T, c=(1,...,1)^T \in \R^d$. Translation by $j=(j_1,...,j_d)$ and dilation by $2^k$ yields $\Psi_{j-2^kc,k}^e(x) := 2^{kd/2} \Psi^e(2^k (x-a) - A (c+2j))$, $2^k \leq j_i \leq 2^{k+1}-1$, $k \in \N_0$, $i=1,...,d$, $e \in E$, that form an orthonormal basis of $L^2([-A,A]^d)$. Then, each $f \in L^2([-A,A]^d)$ has the expansion
\begin{equation}
 f_t =(f_t, \Psi^*) \Psi^* + \sum\limits_{e \in E} \sum\limits_{k=0}^\infty \sum\limits_{\substack{2^k \leq j_i \leq 2^{k+1}-1 \\ i=1,\cdots,d}} (f_t, \Psi_{j-2^kc,k}^e) \Psi_{j-2^kc,k}^e, \label{eq:expansion}
\end{equation}
in the sense of $L^2([-A,A]^d)$ convergence, where
\begin{eqnarray}
 \Psi^*(x) &:=& \frac{1}{(2A)^{d/2}},\quad x \in [-A,A]^d, \\
 (f_t, \Psi_{j-2^kc,k}^e) &:=& \int\limits_{[-A,A]^d} f(x)\Psi_{j-2^kc,k}^e(x) dx.
\end{eqnarray}
It can be shown that any function $f \in L^{\max\{s,p\}}([-A,A]^d)$ for some $p>1$ can be represented by a wavelet series of the form (\ref{eq:expansion}) in the sense of $L^{\max\{s,p\}}([-A,A]^d)$ convergence. However, there exist examples of functions for which (\ref{eq:expansion}) does not hold in particular for $s = p=1$, cf.~\cite{DL92}, p.~7. Therefore, we restrict our setting to kernel functions $f_t \in L^{\max\{s,p\}}([-A,A]^d)$, $p>1$.

As noted, we want to use the expansion~$(\ref{eq:expansion})$ in order to approximate the kernel functions $f_t$ by truncating the (potentially) infinite sum to a finite number of summands. The goal is then to find an upper bound for the approximation error.

\subsubsection{Approximation by cutting off at a certain detail level} \label{subsec:cut_off}
{\ }\newline
Consider a kernel function $f_t \in L^{\max\{s,p\}}([-A,A]^d)$ for some $p>1$ with corresponding Haar series
$$f_t =(f_t, \Psi^*) \Psi^* + \sum\limits_{e \in E} \sum\limits_{k=0}^\infty \sum\limits_{\substack{2^k \leq j_i \leq 2^{k+1}-1 \\ i=1,\cdots,d}} (f_t, \Psi_{j-2^kc,k}^e) \Psi_{j-2^kc,k}^e.$$
The idea is now to cut off this series at a certain detail level $k=n$, that is to approximate the kernel function $f_t$ by
$$\tilde{f}_{t,cut}^{(n)} =(f_t, \Psi^*) \Psi^* + \sum\limits_{e \in E} \sum\limits_{k=0}^n \sum\limits_{\substack{2^k \leq j_i \leq 2^{k+1}-1 \\ i=1,\cdots,d}} (f_t, \Psi_{j-2^kc,k}^e) \Psi_{j-2^kc,k}^e.$$

The following lemma provides an upper bound for the approximation error of bounded kernel functions by applying the cut-off truncation method.

\begin{lemma}\label{lemma:error_bound1}
 Let $s > 0$ and $d > s$. Assume that $M_t:=\sup\limits_{x \in [-A,A]^d} |f_t(x)| < \infty$. Then for $n \in \N_0$
	$$\Vert f_t - \tilde{f}_{t,cut}^{(n)} \Vert_{L^s} \leq
	\begin{cases}
	 \left(\frac{2^d-1}{2^{d-s}-1}\right)^{1/s} \cdot d^{1/s} \cdot M_t \cdot (2A)^{d/s} \cdot \left(\frac{1}{2^{d/s-1}}\right)^{n}, & 0<s< 1, \\
	 \frac{2^d-1}{2^{d/s-1}-1} \cdot d \cdot M_t \cdot (2A)^{d/s} \cdot \left(\frac{1}{2^{d/s-1}}\right)^n, & s \geq 1.
	\end{cases}$$
\end{lemma}

\begin{proof}
 Let $s \geq 1$. We have
\begin{eqnarray}
 \Vert f_t - \tilde{f}_{t,cut}^{(n)} \Vert_{L^s} &=& \left\Vert \sum\limits_{e \in E} \sum\limits_{k=n+1}^\infty \sum\limits_{\substack{2^k \leq j_i \leq 2^{k+1}-1 \\ i=1,\cdots,d}} (f_t, \Psi_{j-2^kc,k}^e) \Psi_{j-2^kc,k}^e \right\Vert_{L^s} \nonumber \\
	&\leq& \sum\limits_{e \in E} \sum\limits_{k=n+1}^\infty \sum\limits_{\substack{2^k \leq j_i \leq 2^{k+1}-1 \\ i=1,\cdots,d}}  |(f_t, \Psi_{j-2^kc,k}^e)| \left\Vert\Psi_{j-2^kc,k}^e \right\Vert_{L^s} \label{eq:greedy}
\end{eqnarray}
Now $|(f_t, \Psi_{j-2^kc,k}^e)|$ can be estimated by
\begin{eqnarray*}
 |(f_t, \Psi_{j-2^kc,k}^e)| &\leq& M_t \int\limits_{[-A,A]^d} |\Psi_{j-2^kc,k}^e(x)|dx
 = M_t \cdot \frac{2^{kd/2}}{(2A)^{d/2}} \cdot 2^{-kd} \cdot (2A)^d \\
 &=& M_t \cdot 2^{-kd/2} \cdot (2A)^{d/2}
\end{eqnarray*}

and $\left\Vert\Psi_{j-2^kc,k}^e \right\Vert_{L^s}$ is equal to
\begin{eqnarray*}
 \left\Vert\Psi_{j-2^kc,k}^e \right\Vert_{L^s} &=& \left(\int_{[-A,A]^d} |\Psi_{j-2^kc,k}^e(x)|^s dx \right)^{1/s}
	= \frac{2^{\frac{kd}{2}}}{(2A)^{d/2}} \cdot 2^{-\frac{kd}{s}} \cdot (2A)^{\frac{d}{s}}.
\end{eqnarray*}

Thus
\begin{eqnarray*}
 \Vert f_t - \tilde{f}_{t,cut}^{(n)} \Vert_{L^s} &\leq& \sum\limits_{e \in E} \sum\limits_{k=n+1}^\infty \sum\limits_{\substack{2^k \leq j_i \leq 2^{k+1}-1 \\ i=1,\cdots,d}}  M_t \cdot 2^{-kd/2} \cdot 2^{\frac{kd}{2}} \cdot 2^{-\frac{kd}{s}}\cdot (2A)^{\frac{d}{s}} \\
 &=& \left( 2^d-1\right) M_t (2A)^{d/s} \sum\limits_{k=n+1}^\infty d\cdot 2^k \cdot 2^{-kd/s} \\
 &=& \left( 2^d-1\right) M_t d (2A)^{d/s}\sum\limits_{k=n+1}^\infty \left( 2^{1-d/s}\right)^k.
\end{eqnarray*}

Since $d > s$, we have $1-d/s < 0$ and by using the geometric series formula, we get
\begin{equation}
 \Vert f_t - \tilde{f}_{t,cut}^{(n)} \Vert_{L^s} \leq \frac{2^d-1}{2^{d/s-1} -1} d M_t (2A)^{d/s}\left(\frac{1}{2^{d/s-1}} \right)^n.
\end{equation}

Now let $0<s < 1$. By Lemma \ref{lemma:inequality}, we have
\begin{eqnarray*}
 \Vert f_t - \tilde{f}_{t,cut}^{(n)} \Vert_{L^s}^s &=& \left\Vert \sum\limits_{e \in E} \sum\limits_{k=n+1}^\infty \sum\limits_{\substack{2^k \leq j_i \leq 2^{k+1}-1 \\ i=1,\cdots,d}} (f_t, \Psi_{j-2^kc,k}^e) \Psi_{j-2^kc,k}^e \right\Vert_{L^s}^s \nonumber \\
	&\leq& \sum\limits_{e \in E} \sum\limits_{k=n+1}^\infty \sum\limits_{\substack{2^k \leq j_i \leq 2^{k+1}-1 \\ i=1,\cdots,d}}  |(f_t, \Psi_{j-2^kc,k}^e)|^s \left\Vert\Psi_{j-2^kc,k}^e \right\Vert_{L^s}^s. \nonumber
\end{eqnarray*}
By using the estimates for the wavelet coefficients and the $L^s$-norms of the wavelets from above, we get
\begin{eqnarray*}
 \Vert f_t - \tilde{f}_{t,cut}^{(n)} \Vert_{L^s}^s &\leq& \left( 2^d-1\right) d \left(M_t (2A)^{d/s}\right)^s \sum\limits_{k=n+1}^\infty \left( 2^{s-d}\right)^k
\end{eqnarray*}
and finally
$$\Vert f_t - \tilde{f}_{t,cut}^{(n)} \Vert_{L^s} \leq \frac{\left(2^d-1\right)^{1/s}}{\left(2^{d-s}-1\right)^{1/s}} \cdot d^{1/s} \cdot M_t \cdot (2A)^{d/s} \cdot \left(\frac{1}{2^{d/s-1}}\right)^{n}.$$\qed
\end{proof}

If we make further assumptions about the kernel function $f_t$, we can improve the rate of convergence of the upper bound.

\begin{corollary}
 Assume that $f_t$ is Hölder-continuous with parameters $C_t$ and ${\gamma_t}$ for all $t \in [-T,T]^q$. Then for $n \in \N_0$
$$\Vert f_t - \tilde{f}_{t,cut}^{(n)} \Vert_{L^s} \leq
	\begin{cases}
	 \frac{1}{2}\left(\frac{2^d-1}{2^{d+s\gamma_t}-1}\right)^{1/s}\cdot d^{1/s+\gamma_t/(2s)} \cdot C_t \cdot (2A)^{d/s+\gamma_t} \cdot \left(\frac{1}{2^{d/s+\gamma_t}}\right)^{n}, & 0<s < 1, \\
	 \frac{2^d-1}{2^{d/s+\gamma_t+1}-2} \cdot d^{1+\gamma_t/2} \cdot C_t \cdot (2A)^{d/s+\gamma_t} \cdot \left(\frac{1}{2^{d/s+\gamma_t}}\right)^{n}, & s \geq 1.
	\end{cases}$$
\end{corollary}

\begin{proof}
 We estimate $|(f_t, \Psi_{j-2^kc,k}^e)|$ as in the preceeding lemma. Let 
 \begin{eqnarray*}
  B &:=& \{x \in [-A,A]^d: \Psi_{j-2^kc,k}^e(x) > 0\}, \\
  C &:=& \{x \in [-A,A]^d: \Psi_{j-2^kc,k}^e(x) < 0\}.
 \end{eqnarray*}
 Then we have
 \begin{eqnarray*}
  &&|(f_t, \Psi_{j-2^kc,k}^e)| \\
 &=& \left\vert \int_{[-A,A]^d} f_t(x) \Psi_{j-2^kc,k}^e(x) dx \right\vert
	= \left\vert \int_B f_t(x) \frac{2^{kd/2}}{(2A)^{d/2}} dx - \int_C f_t(x) \frac{2^{kd/2}}{(2A)^{d/2}} dx \right\vert \\
	&=& \frac{2^{kd/2}}{(2A)^{d/2}} \left\vert \int_B f_t(x) dx - \int_C f_t(x) dx \right\vert \\
	&\leq& \frac{2^{kd/2}}{(2A)^{d/2}} \max \left\{ \left\vert \int_B f_t(x) dx - \max\limits_{x \in C}(f_t(x)) \vert C \vert \right\vert, \left\vert \int_B f_t(x) dx - \min\limits_{x \in C}(f_t(x)) \vert C \vert \right\vert \right\}.
 \end{eqnarray*}
Here, $\vert C \vert$ denotes the volume of $C$. We now estimate the two quantities in the maximum. We have
\begin{eqnarray*}
 &&\left\vert \int_B f_t(x) dx - \max\limits_{x \in C}(f_t(x)) \vert C \vert \right\vert \\
 &\leq& \max \left\{ \left\vert \max\limits_{x \in B}(f_t(x)) \vert B \vert - \max\limits_{x \in C}(f_t(x)) \vert C \vert \right\vert , \left\vert \min\limits_{x \in B}(f_t(x)) \vert B \vert - \max\limits_{x \in C}(f_t(x)) \vert C \vert \right\vert \right\}
\end{eqnarray*}
and
\begin{eqnarray*}
 &&\left\vert \int_B f_t(x) dx - \min\limits_{x \in C}(f_t(x)) \vert C \vert \right\vert \\
 &\leq& \max \left\{ \left\vert \max\limits_{x \in B}(f_t(x)) \vert B \vert - \min\limits_{x \in C}(f_t(x)) \vert C \vert \right\vert , \left\vert \min\limits_{x \in B}(f_t(x)) \vert B \vert - \min\limits_{x \in C}(f_t(x)) \vert C \vert \right\vert \right\}.
\end{eqnarray*}
Therefore we get
\begin{eqnarray*}
 &&|(f_t, \Psi_{j-2^kc,k}^e)| \\
 &\leq& \frac{2^{kd/2}}{(2A)^{d/2}} \cdot \\
 && \hspace*{0.5cm} \cdot \max\left\{ \left\vert \max\limits_{x \in B}(f_t(x)) \vert B \vert - \max\limits_{x \in C}(f_t(x)) \vert C \vert \right\vert , \left\vert \min\limits_{x \in B}(f_t(x)) \vert B \vert - \max\limits_{x \in C}(f_t(x)) \vert C \vert \right\vert \right., \\
 && \hspace*{1.7cm}\left. \left\vert \max\limits_{x \in B}(f_t(x)) \vert B \vert - \min\limits_{x \in C}(f_t(x)) \vert C \vert \right\vert , \left\vert \min\limits_{x \in B}(f_t(x)) \vert B \vert - \min\limits_{x \in C}(f_t(x)) \vert C \vert \right\vert \right\}.
\end{eqnarray*}
Let $x_1 \in B$ and $x_2 \in C$ such that the maximum in the last inequality is attained. Furthermore, we have
\begin{equation*}
 \vert B \vert = \vert C \vert = \frac{1}{2} \cdot (2A)^d \cdot 2^{-kd}.
\end{equation*}

Then, since $f_t$ is H\"older-continuous, we get
\begin{eqnarray*}
 &&|(f_t, \Psi_{j-2^kc,k}^e)| \\
&\leq& \frac{2^{kd/2}}{(2A)^{d/2}} \left\vert f_t(x_1) \cdot \vert B \vert - f_t(x_2) \cdot \vert C \vert \right\vert 
	= \frac{2^{kd/2}}{(2A)^{d/2}} \cdot \frac{1}{2} (2A)^d \cdot 2^{-kd} \vert f_t(x_1) - f_t(x_2)\vert \\
	&\leq& \frac{2^{kd/2}}{(2A)^{d/2}} \cdot \frac{1}{2} (2A)^d \cdot 2^{-kd} C_t \Vert x_1-x_2 \Vert^{\gamma_t} 
	\leq \frac{2^{kd/2}}{(2A)^{d/2}} \cdot \frac{1}{2} (2A)^d \cdot 2^{-kd} C_t \left(2A\sqrt{d} 2^{-k}\right)^{\gamma_t} \\
	&=& A^{d/2+\gamma_t} C_t 2^{d/2+\gamma_t-1} d^{\gamma_t/2} \left(\frac{1}{2^{d/2+\gamma_t}}\right)^k.
\end{eqnarray*}
The remainder of the proof is analogous to the one of Lemma \ref{lemma:error_bound1}.\qed
\end{proof}

\begin{corollary}
 Assume that $f_t$ is differentiable with $||\nabla f_t(x)||_2 \leq C_t$ for all \linebreak $x \in [-A,A]^d$, $C_t > 0$ and $t \in [-T,T]^q$. Then for $n \in \N_0$
$$\Vert f_t - \tilde{f}_{t,cut}^{(n)} \Vert_{L^s} \leq
	\begin{cases}
	 \frac{1}{2}\left(\frac{2^d-1}{2^{d+s}-1}\right)^{1/s}\cdot d^{3/(2s)} \cdot C_t \cdot (2A)^{d/s+1} \cdot \left(\frac{1}{2^{d/s+1}}\right)^{n}, & 0<s < 1, \\
	 \frac{2^d-1}{2^{d/s+2}-2} \cdot d^{3/2} \cdot C_t \cdot (2A)^{d/s+1} \cdot \left(\frac{1}{2^{d/s+1}}\right)^{n}, & s \geq 1.
	\end{cases}$$
\end{corollary}

\subsubsection{Near best $n$-term approximation}
{\ }\newline
Taking a wavelet basis for $L^{\max\{s,p\}}([-A,A]^d)$, $p>1$, has advantages in particular in the representation of functions with discontinuities and sharp peaks, that is functions with a certain local behavior. By simply cutting of at a certain detail level, this advantage is not honored. In view of (\ref{eq:greedy}), we may expect that we have to calculate less Haar coefficients if we approximate the kernel function $f_t$ by a truncated Haar series $\tilde{f}_t^{(n)}$ that contains those $n$ summands $(f, \Psi_{j-2^kc,k}^e) \Psi_{j-2^kc,k}^e$ with the largest values $\Vert (f, \Psi_{j-2^kc,k}^e) \Psi_{j-2^kc,k}^e \Vert_{L^s}=|(f, \Psi_{j-2^kc,k}^e)| \left\Vert\Psi_{j-2^kc,k}^e \right\Vert_{L^s}$. An approach to use such a truncation in order to approximate functions is presented in \cite{DeV98}, p.~114 ff. We summarize the main statements.

Consider a function $S$ defined by
\begin{equation}
S = \sum\limits_{(e,k,j) \in \Xi} a_{j,k}^e \Psi_{j-2^kc,k}^e, \quad a_{j,k}^e \in \R, \quad \forall (e,j,k) \in \Xi \label{eq:S}
\end{equation}
where $\Xi := \{(e,k,j): e \in E, k \in \N_0, 2^k \leq j_i \leq 2^{k+1}-1, i=1,...,d\}$ and $\# \Xi \leq n$ for some $n \in \N$. Hence, $(\ref{eq:S})$ is a linear combination of $n$ Haar wavelets. We denote by $\Sigma_n$ the set of all the functions $S$ defined as in $(\ref{eq:S})$ and let
$$\sigma_n^s(f_t) := \inf\limits_{S \in \Sigma_n} ||f - S||_{L^s}.$$

Now, we truncate the wavelet expansion $(\ref{eq:expansion})$ of $f_t$ by taking those $n$ summands for which the absolute value of $|(f, \Psi_{j-2^kc,k}^e)| \left\Vert\Psi_{j-2^kc,k}^e \right\Vert_{L^s}$ is largest and denote the truncated sum by $\tilde{f}_t^{(n)}$. In \cite{Tem98}, the following theorem was proven which shows that this truncation is a near best $n$-term approximation.

\begin{theorem} \label{th:temlyakov}
 Let $1 < s < \infty$. Then for any $f \in L^s([-A,A]^d)$ we have
 $$\Vert f - \tilde{f}^{(n)} \Vert_{L^s} \leq C_1(s,d,A) \sigma_n^s(f_t)$$
 with a constant $C_1(p,d,A) \geq 0$ only depending on $s$, $d$ and $A$.
\end{theorem}

If the sequence $\left\{|(f, \Psi_{j-2^kc,k}^e)| \left\Vert\Psi_{j-2^kc,k}^e \right\Vert_{L^s}\right\}_{j,k}$ is in the Lorentz space $wl_\tau$, \linebreak $0<\tau<\infty$, that is
\begin{equation}
 \# \{(j,k,e): \{\Vert(f, \Psi_{j,k}^e)\Psi_{j,k}^e \Vert_{L^s} > \varepsilon\} \leq \left(\frac{M}{\varepsilon}\right)^\tau, \quad \forall \varepsilon > 0, \quad M \geq 0, \label{eq:lorentz}
\end{equation}
with a certain additional condition on $\tau$, then $\sigma_n^s(f_t)$ can be bounded from above as shown in \cite{DeV98}, p.~116.

\begin{theorem} \label{th:devore}
 Let $1 < s < \infty$ and $f \in L^s([-A,A]^d)$. Furthermore, let 
$$\left\{\left\Vert(f, \Psi_{j-2^kc,k}^e)\Psi_{j-2^kc,k}^e \right\Vert_{L^s}\right\}_{j,k} \in wl_\tau$$
and $u > 0$ with $1/\tau = u + 1/s$. Then
$$\sigma_n^s(f) \leq C_2(s,d,A) M \left(\frac{1}{n}\right)^{u}, \quad n = 1,2,\cdots,$$
with a constant $C_2(s,d,A) \geq 0$ only depending on $s$, $d$ and $A$ and $M$ being a constant satisfying $(\ref{eq:lorentz})$.
\end{theorem}

Combining Theorem \ref{th:temlyakov} and Theorem \ref{th:devore} yields an upper bound of the approximation error by using the near best $n$-term approximation with a rate of convergence of $O\left((1/n)^s\right)$. In \cite{KSS09a}, we obtained the following formulas for $C_1(s,d,A)$ and $C_2(s,d,A)$ in the case that $f \in L^s([0,1]^d)$:
\begin{eqnarray*}
 C_1(s,d) &=& \left(2+\left(1-2^{-\frac{d}{s}}\right)^{-2}\right) \left(\left(2^d-1\right) \left(\max\left(s,\frac{s}{s-1}\right) -1 \right)\right)^2, \\
 C_2(s,d) &=& \frac{2 \left(2^{\tau/s}-1\right)}{\left(1-\left(\frac{1}{2}\right)^{d/s}\right)\left(1-\left(\frac{1}{2}\right)^\tau \right)^{1/s} \left(1-2^{\tau/s-1}\right)}.
\end{eqnarray*}

We note that on the one hand, these constants are not sharp and may be quite large, and on the other hand, we would have to find a value of $M$ in (\ref{eq:lorentz}) as small as possible for each kernel function $f_t$ or for certain classes of kernel functions, which is not so easy to determine. Therefore, we suggest an approach which is not based on the error estimate with those constants, but still determines an approximation at least close to the near best $n$-term approximation while keeping the desired level of accuracy.

It is clear that $\Vert (f, \Psi_{j-2^kc,k}^e) \Psi_{j-2^kc,k}^e \Vert_{L^s}$ goes to zero as the detail level $k$ goes to infinity. This means that the $n$ largest values $\Vert (f, \Psi_{j-2^kc,k}^e) \Psi_{j-2^kc,k}^e \Vert_{L^s}$ are likely to be found for small values of $k$.

We now assume that $s \geq 1$, $d > s$ and $M_t:=\sup\limits_{x \in [-A,A]^d} |f_t(x)| < \infty$ and choose $\varepsilon > 0$ as the desired level of accuracy. The following derivations are analogous for the case $0 < s < 1$.

From Lemma \ref{lemma:error_bound1} and its proof, we get
\begin{eqnarray*}
 \Vert f_t - \tilde{f}_{t,cut}^{(n)} \Vert_{L^s} 
 &\leq& \sum\limits_{e \in E} \sum\limits_{k=n+1}^\infty \sum\limits_{\substack{2^k \leq j_i \leq 2^{k+1}-1 \\ i=1,\cdots,d}}  \Vert (f, \Psi_{j-2^kc,k}^e) \Psi_{j-2^kc,k}^e \Vert_{L^s} \\
 &\leq& \frac{2^d-1}{2^{d/s-1}-1} \cdot d \cdot M_t \cdot (2A)^{d/s} \cdot \left(\frac{1}{2^{d/s-1}}\right)^n
 \stackrel{!}{\leq} \varepsilon
\end{eqnarray*}

if and only if

$$ n \geq \frac{\ln(\varepsilon (2^{d/s-1}-1)) - \ln ((2^d-1)dM_t(2A)^{d/s})}{-\ln(2^{d/s-1})}. $$

We take

$$ m_t: = \left\lceil \frac{\ln(\varepsilon (2^{d/s-1}-1)) - \ln ((2^d-1)dM_t(2A)^{d/s})}{-\ln(2^{d/s-1})} \right\rceil,$$
where $\lceil x \rceil$ is the integral part of $x$, as our minimal detail level to obtain the desired level of accuracy and add $l \in \N_0$ detail levels to the truncated wavelet series at detail level $m_t$, that is we consider
\begin{equation}
 \tilde{f}_{t,cut}^{(m_t+l)} = (f_t, \Psi^*) \Psi^* + \sum\limits_{e \in E} \sum\limits_{k=0}^{m_t+l} \sum\limits_{\substack{2^k \leq j_i \leq 2^{k+1}-1 \\ i=1,\cdots,d}} (f_t, \Psi_{j-2^kc,k}^e) \Psi_{j-2^kc,k}^e \label{eq:m+l}.
\end{equation}
Furthermore, we define
\begin{equation}
 C := \Vert (f_t, \Psi^*) \Psi^* \Vert_{L^s}+ \sum\limits_{e \in E} \sum\limits_{k=0}^{m_t+l} \sum\limits_{\substack{2^k \leq j_i \leq 2^{k+1}-1 \\ i=1,\cdots,d}} \Vert (f_t, \Psi_{j-2^kc,k}^e) \Psi_{j-2^kc,k}^e \Vert_{L^s} \label{eq:number}
\end{equation}
and
$$ D := \Vert (f_t, \Psi^*) \Psi^* \Vert_{L^s}+ \sum\limits_{e \in E} \sum\limits_{k=m_t+l}^{\infty} \sum\limits_{\substack{2^k \leq j_i \leq 2^{k+1}-1 \\ i=1,\cdots,d}} \Vert (f_t, \Psi_{j-2^kc,k}^e) \Psi_{j-2^kc,k}^e \Vert_{L^s}.$$
By Lemma \ref{lemma:error_bound1}, the corresponding level of accuracy of (\ref{eq:m+l}) is at least
$$\varepsilon_t^* := \frac{2^d-1}{2^{d/s-1}-1} \cdot d \cdot M_t \cdot (2A)^{d/s} \cdot \left(\frac{1}{2^{d/s-1}}\right)^{m_t+l}.$$

Now we take the $n$ largest summands from (\ref{eq:number}), that is from $$\{\Vert (f_t, \Psi^*) \Psi^* \Vert_{L^s}\} \bigcup \{\Vert (f_t, \Psi_{j-2^kc,k}^e) \Psi_{j-2^kc,k}^e \Vert_{L^s}\}_{j,k},$$ and denote them by $a_1$,...,$a_n$. The remaining summands are denoted by $a_{n+1}$, $a_{n+2}$, ..., and the corresponding summands from (\ref{eq:m+l}) by $b_1$,...,$b_n$, $b_{n+1}$, $b_{n+2}$, ... . The number $n \in \N$ is chosen to be the smallest number such that
$$ C - \sum\limits_{i=1}^n a_i \leq \varepsilon - \varepsilon_t^*.$$

We define
$$\tilde{f}_t^{(n)} := \sum\limits_{i=1}^n b_i$$
which is close to the near best $n$-term approximation if $l$ is chosen large enough since 
$$\Vert (f_t,\Psi_{j-2^{m_t+l},m_t+l}^e) \Psi_{j-2^{m_t+l},m_t+l}^e \Vert_{L^s}$$ goes to zero as $l$ goes to infinity.

Then we have
$$\Vert f - \tilde{f}_t^{(n)} \Vert_{L^s} = \left\Vert \sum\limits_{i=n+1}^\infty b_i \right\Vert_{L^s} \leq \sum\limits_{i=n+1}^\infty a_i = C - \sum\limits_{i=1}^n a_i + D \leq \varepsilon - \varepsilon^* + \varepsilon^* = \varepsilon. $$

\subsubsection{Implementation}
{\ }\newline
When implementing the wavelet approach, one more problem has to be considered which we discuss now.

Let $I$ be the set of the indices $(e,j,k)$ for which the summands $(f_t, \Psi_{j-2^kc,k}^e) \Psi_{j-2^kc,k}^e$ are part of the approximation $\tilde{f}_t^{(n)}$. Then we can write
$$ \tilde{f}_t^{(n)} = (f_t, \Psi^*) \Psi^* + \sum\limits_{(e,k,j) \in I} (f_t, \Psi_{j-2^kc,k}^e) \Psi_{j-2^kc,k}^e $$
if $(f_t, \Psi^*) \Psi^*$ is included in the truncated series or
$$ \tilde{f}_t^{(n)} = \sum\limits_{(e,k,j) \in I} (f_t, \Psi_{j-2^kc,k}^e) \Psi_{j-2^kc,k}^e $$
if it is not included.

In order to approximate the random field $X$, we use
$$\tilde{X}^{(n)}(t) = (f_t, \Psi^*) \cdot \frac{\Lambda([-A,A]^d)}{(2A)^{d/2}} + \sum\limits_{(e,k,j) \in I} (f_t, \Psi_{j-2^kc,k}^e) \int\limits_{[-A,A]^d} \Psi_{j-2^kc,k}^e \Lambda(dx) $$
or
$$\tilde{X}^{(n)}(t) = \sum\limits_{(e,k,j) \in I} (f_t, \Psi_{j-2^kc,k}^e) \int\limits_{[-A,A]^d} \Psi_{j-2^kc,k}^e \Lambda(dx) $$
if, again, $(f_t, \Psi^*) \Psi^*$ is not included in the truncated series.

Since the Haar wavelets $\Psi_{j-2^kc,k}^e$ are simple step functions, the integrals 
$$\int\limits_{[-A,A]^d} \Psi_{j-2^kc,k}^e \Lambda(dx)$$ can be easily simulated although they are not independent: Let $z \in \N$ be the finest detail level of the wavelet approximation. Then all of these integrals can be built up from
$$\int\limits_{\left[-A+k\frac{A}{2^z},-A+(k+1)\frac{A}{2^z}\right)^d} \Lambda(dx) = \Lambda\left(\left[-A+k\frac{A}{2^z},-A+(k+1)\frac{A}{2^z}\right)^d\right), \hspace*{0.1cm} k=0,...,2^{z+1}-1$$
which are independent because the sets $\left[-A+k\frac{A}{2^z},-A+(k+1)\frac{A}{2^z}\right)^d$ are disjoint.

However, the wavelet coefficients $(f_t, \Psi_{j-2^kc,k}^e)$ cause problems if no closed formula of the integral of the kernel functions $f_t$ over cubes is known. In this case, they have to be determined numerically by using the fast wavelet transform (see for instance \cite{Urb08}, pp.~134).

This results in a further approximation error which we need to estimate. When the detail level at which the wavelet series is cut off is equal to $n$, the input vector of the fast wavelet transform consists of integrals of the form
$$\int_{cube} \frac{2^{(n+1)d/2}}{(2A)^{d/2}} f_t(x) dx,$$
where $cube$ is a cube of side length $2^{-d(n+1)}$.

We now assume that we have calculated these integrals with a precision of $\delta > 0$ and denote the wavelet coefficients computed by the fast wavelet transform by $\widehat{(f_t, \Psi^*)}$ and $\widehat{(f_t, \Psi_{j-2^kc,k}^e)}$.

When applying the fast wavelet transform, the value of each integral is used $2^d-1$ times at each detail level $k$ to calculate the wavelet coefficients $(f_t, \Psi_{j-2^kc,k}^e)$, $e \in E$. There are $2^{(n+1)d}$ such integrals.

When $s \geq 1$, the precision of
$$\widehat{\tilde{f}_{t,cut}^{(n)}} =\widehat{(f_t, \Psi^*)} \Psi^* + \sum\limits_{e \in E} \sum\limits_{k=0}^n \sum\limits_{\substack{2^k \leq j_i \leq 2^{k+1}-1 \\ i=1,\cdots,d}} \widehat{(f_t, \Psi_{j-2^kc,k}^e)} \Psi_{j-2^kc,k}^e$$
to approximate
\begin{equation}
 \tilde{f}_{t,cut}^{(n)} =(f_t, \Psi^*) \Psi^* + \sum\limits_{e \in E} \sum\limits_{k=0}^n \sum\limits_{\substack{2^k \leq j_i \leq 2^{k+1}-1 \\ i=1,\cdots,d}} (f_t, \Psi_{j-2^kc,k}^e) \Psi_{j-2^kc,k}^e \label{eq:rf_wavelet}
\end{equation}
is
\begin{eqnarray}
&& \hspace*{-1.2cm}\Vert \widehat{\tilde{f}_{t,cut}^{(n)}} - \tilde{f}_{t,cut}^{(n)} \Vert_{L^s}
\leq |\widehat{(f_t, \Psi^*)}-(f_t, \Psi^*)| \Vert \Psi^* \Vert_{L^s} \nonumber \\
 &&\hspace*{2.3cm}+ \sum\limits_{e \in E} \sum\limits_{k=0}^n \sum\limits_{\substack{2^k \leq j_i \leq 2^{k+1}-1 \\ i=1,\cdots,d}} |\widehat{(f_t, \Psi_{j-2^kc,k}^e)}-(f_t, \Psi_{j-2^kc,k}^e)| \Vert \Psi_{j-2^kc,k}^e \Vert_{L^s} \nonumber\\
&\leq& 2^{(n+1)d} \frac{(2A)^{d/s-d/2}}{2^{(n+1)d/2}} \delta + (2^d-1) \sum\limits_{k=0}^n 2^{(n+1)d} \frac{2^{kd/2}}{2^{(n+1)d/2}} \delta \cdot (2A)^{d/s-d/2} \cdot 2^{kd/2-kd/s} \nonumber \\
&\leq& \begin{cases}
        (2A)^{d/s-d/2} 2^{(n+1)d/2}  \left(1+(2^d-1) \frac{2^{(d-d/s)(n+1)}-1}{2^{d-d/s}-1}\right)\delta, & s > 1,\\
	(2A)^{d/2} 2^{(n+1)d/2} \left(1+(2^d-1)(n+1)\right)\delta, & s = 1.
       \end{cases} \label{eq:delta}
\end{eqnarray}

When $0 < s < 1$, we can use Lemma \ref{lemma:inequality} and get
\begin{eqnarray*}
 \Vert \widehat{\tilde{f}_{t,cut}^{(n)}} - \tilde{f}_{t,cut}^{(n)} \Vert_{L^s} \nonumber
&\leq& \left( |\widehat{(f_t, \Psi^*)}-(f_t, \Psi^*)|^s \Vert \Psi^* \Vert_{L^s}^s \right. \nonumber\\
&& \hspace*{-0.7cm}+ \sum\limits_{e \in E} \sum\limits_{k=0}^n \sum\limits_{\substack{2^k \leq j_i \leq 2^{k+1}-1 \\ i=1,\cdots,d}}\left. |\widehat{(f_t, \Psi_{j-2^kc,k}^e)}-(f_t, \Psi_{j-2^kc,k}^e)|^s \Vert \Psi_{j-2^kc,k}^e \Vert_{L^s}^s \right)^{\frac{1}{s}} \\
& \leq& (2A)^{d/s-d/2} 2^{(n+1)d/2}  \left(1+(2^d-1) \frac{2^{(ds-d)(n+1)}-1}{2^{ds-d}-1}\right)^{1/s}\delta.
\end{eqnarray*}

Let $\varepsilon > 0$. We choose $\varepsilon_1 > 0$ and $\varepsilon_2 > 0$ such that $\varepsilon_1+\varepsilon_2 = \varepsilon$ if $s \geq 1$ and $\varepsilon_1^s+\varepsilon_2^s = \varepsilon^s$ if $0 < s < 1$. Furthermore, we choose the detail level $n$ so large (by using the formulas in Section \ref{subsec:cut_off}) such that
$$ \Vert f_t - \tilde{f}_{t,cut}^{(n)} \Vert_{L^s} \leq \varepsilon_1.$$
We approximate the elements of the input vector for the fast wavelet transform with a precision of
$$\delta = \begin{cases}
	    \frac{\varepsilon_2}{(2A)^{d/s-d/2} 2^{(n+1)d/2}  \left(1+(2^d-1) \frac{2^{(ds-d)(n+1)}-1}{2^{ds-d}-1}\right)^{1/s}}, & 0 < s < 1, \\
	    \frac{\varepsilon_2}{(2A)^{d/2} 2^{(n+1)d/2} \left(1+(2^d-1)(n+1)\right)}, & s=1, \\
            \frac{\varepsilon_2}{(2A)^{d/s-d/2} 2^{(n+1)d/2}  \left(1+(2^d-1) \frac{2^{(d-d/s)(n+1)}-1}{2^{d-d/s}-1}\right)}, & s > 1.
           \end{cases}$$
Then we have for $s \geq 1$
$$\Vert f_t - \widehat{\tilde{f}_{t,cut}^{(n)}}\Vert_{L^s} \leq \Vert f_t - \tilde{f}_{t,cut}^{(n)} \Vert_{L^s} + \Vert \tilde{f}_{t,cut}^{(n)}- \widehat{\tilde{f}_{t,cut}^{(n)}} \Vert_{L^s} \leq \varepsilon_1+\varepsilon_2 = \varepsilon,$$
and for $0<s < 1$
$$\Vert f_t - \widehat{\tilde{f}_{t,cut}^{(n)}}\Vert_{L^s} \leq \left(\Vert f_t - \tilde{f}_{t,cut}^{(n)} \Vert_{L^s}^s + \Vert \tilde{f}_{t,cut}^{(n)}- \widehat{\tilde{f}_{t,cut}^{(n)}} \Vert_{L^s}^s\right)^{1/s} = \varepsilon. $$

We summarize this result in the following algorithm.

\vspace*{0.5cm}
\underline{\textbf{Algorithm}}

Let $M_t:=\sup\limits_{x \in [-A,A]^d} |f_t(x)| < \infty$ and $d > s$. Choose $\varepsilon > 0$ as the desired level of accuracy. Choose $\varepsilon_1, \varepsilon_2>0$ such that $\varepsilon = \varepsilon_1 + \varepsilon_2$ if $s \geq 1$ and $\varepsilon = (\varepsilon_1^s+\varepsilon_2^s)^{1/s}$ if $0<s<1$.
\begin{enumerate}
 \item Let 
	$$m_t: = \begin{cases}
	 \left\lceil \frac{\ln(\varepsilon_1 (2^{d-s}-1)^{1/s}) - \ln ((2^d-1)^{1/s} d^{1/s} M_t(2A)^{\frac{d}{s}})}{(1-d/s)\ln(2)} \right\rceil, & 0<s< 1, \\
	 \left\lceil \frac{\ln(\varepsilon_1(2^{d/s-1}-1)) - \ln ((2^d-1)dM_t(2A)^{d/s})}{(1-d/s)\ln(2)} \right\rceil, & s \geq 1.
	\end{cases}$$
	and choose a number $l \in \N_0$ that increases the detail level $m_t$.
 \item Calculate the wavelet coefficients for
	$$ (f_t, \Psi^*) \Psi^* + \sum\limits_{e \in E} \sum\limits_{k=0}^{m_t+l} \sum\limits_{\substack{2^k \leq j_i \leq 2^{k+1}-1 \\ i=1,\cdots,d}}  (f_t, \Psi_{j-2^kc,k}^e) \Psi_{j-2^kc,k}^e $$ using the fast wavelet transform with a precision of
	$$\delta = \begin{cases}
	    \frac{\varepsilon_2}{(2A)^{d/s-d/2} 2^{(n+1)d/2}  \left(1+(2^d-1) \frac{2^{(ds-d)(n+1)}-1}{2^{ds-d}-1}\right)^{1/s}}, & 0 < s < 1, \\
	    \frac{\varepsilon_2}{(2A)^{d} 2^{(n+1)d/2}  \left(1+(2^d-1) \frac{2^{-d/2\cdot(n+1)}-1}{2^{-d/2}-1}\right)}, & s=1, \\
            \frac{\varepsilon_2}{(2A)^{d/s-d/2} 2^{(n+1)d/2}  \left(1+(2^d-1) \frac{2^{(d-d/s)(n+1)}-1}{2^{d-d/s}-1}\right)}, & s > 1.
           \end{cases}$$
 \item Take the $n$ largest summands from
$$C := \Vert \widehat{(f_t, \Psi^*)} \Psi^* \Vert_{L^s}+ \sum\limits_{e \in E} \sum\limits_{k=0}^{m_t+l} \sum\limits_{\substack{2^k \leq j_i \leq 2^{k+1}-1 \\ i=1,\cdots,d}} \Vert \widehat{(f_t, \Psi_{j-2^kc,k}^e)} \Psi_{j-2^kc,k}^e \Vert_{L^s}$$
 and denote them by $a_1$,...,$a_n$. The corresponding summands from
 $$\widehat{(f_t, \Psi^*)} \Psi^* + \sum\limits_{e \in E} \sum\limits_{k=0}^{m_t+l} \sum\limits_{\substack{2^k \leq j_i \leq 2^{k+1}-1 \\ i=1,\cdots,d}} \widehat{(f_t, \Psi_{j-2^kc,k}^e)} \Psi_{j-2^kc,k}^e$$
 are denoted by $b_1$,...,$b_n$. Choose the number $n$ such that
 $$ C - \sum\limits_{i=1}^n a_i \leq \varepsilon_1 - \varepsilon_t^*,$$
 where $$\varepsilon_t^* = \begin{cases}
	 \left(\frac{2^d-1}{2^{d-s}-1}\right)^{1/s} \cdot d^{1/s} \cdot M_t \cdot (2A)^{d/s} \cdot \left(\frac{1}{2^{d/s-1}}\right)^{m_t+l}, & 0<s< 1, \\
	 \frac{2^d-1}{2^{d/s-1}-1} \cdot d \cdot M_t \cdot (2A)^{d/s} \cdot \left(\frac{1}{2^{d/s-1}}\right)^{m_t+l}, & s \geq 1.
	\end{cases}$$.
\item Take $\tilde{f}_t^{(n)} = \sum_{i=1}^n b_i$ as the approximation for $f_t$.
\end{enumerate}

\begin{remark}
 \begin{enumerate}
  \item Assume that $d \in \N$ and $f_t$ is Hölder-continuous with parameters $C_t$ and ${\gamma_t}$ for $t \in [-T,T]^q$. Then the algorithm can be applied with $m_t$ and $\varepsilon_t^*$ replaced by
	\begin{eqnarray*}
	 m_t &:=& \begin{cases}
	        \left\lceil \frac{\ln(2\varepsilon_1 (2^{d+s\gamma_t}-1)^{1/s}) - \ln ((2^d-1)^{1/s}d^{1/s+\gamma_t/(2s)}C_t(2A)^{d/s+\gamma_t})}{-\ln(2^{d/s+\gamma_t})} \right\rceil, & 0 < s < 1,\\
		\left\lceil \frac{\ln(\varepsilon_1 (2^{d/s+\gamma_t+1}-2)) - \ln ((2^d-1)d^{1+\gamma_t/2}C_t(2A)^{d/s+\gamma_t})}{-\ln(2^{d/s+\gamma_t})} \right\rceil, & s \geq 1,
	       \end{cases} \\
	\varepsilon_t^* &:=& \begin{cases}
	 \frac{1}{2}\left(\frac{2^d-1}{2^{d+s\gamma_t}-1}\right)^{1/s}\cdot d^{\frac{1}{s}+\frac{\gamma_t}{2s}} \cdot C_t \cdot (2A)^{\frac{d}{s}+\gamma_t} \cdot \left(\frac{1}{2^{d/s+\gamma_t}}\right)^{m_t+l}, & 0<s < 1, \\
	 \frac{2^d-1}{2^{d/s+\gamma_t+1}-2} \cdot d^{1+\gamma_t/2} \cdot C_t \cdot (2A)^{d/s+\gamma_t} \cdot \left(\frac{1}{2^{d/s+\gamma_t}}\right)^{m_t+l}, & s \geq 1.
	\end{cases}
	\end{eqnarray*}
  \item Assume that $d \in \N$ and $f_t$ is differentiable with $||\nabla f_t(x)||_2 \leq C_t$ for all $x \in \text{supp}(f_t)$ with $C_t > 0$ and $t \in [-T,T]^q$. Then the algorithm can be applied with $m_t$ and $\varepsilon_t^*$ replaced by
	\begin{eqnarray*}
	 m_t &:=& \begin{cases}
	        \left\lceil \frac{\ln(2\varepsilon_1 (2^{d+s}-1)^{1/s}) - \ln ((2^d-1)^{1/s}d^{3/(2s)}C_t(2A)^{d/s+1})}{-\ln(2^{d/s+1})} \right\rceil, & 0 < s < 1,\\
		\left\lceil \frac{\ln(\varepsilon_1 (2^{d/s+2}-2)) - \ln ((2^d-1)d^{3/2}C_t(2A)^{d/s+1})}{-\ln(2^{d/s+1})} \right\rceil, & s \geq 1,
	       \end{cases} \\
	\varepsilon_t^* &:=& \begin{cases}
	 \frac{1}{2}\left(\frac{2^d-1}{2^{d+s}-1}\right)^{1/s}\cdot d^{1/s+1/(2s)} \cdot C_t \cdot (2A)^{d/s+1} \cdot \left(\frac{1}{2^{d/s+1}}\right)^{m_t+l}, & 0<s < 1, \\
	 \frac{2^d-1}{2^{d/s+2}-2} \cdot d^{3/2} \cdot C_t \cdot (2A)^{d/s+1} \cdot \left(\frac{1}{2^{d/s+1}}\right)^{m_t+l}, & s \geq 1.
	\end{cases}
	\end{eqnarray*}
 \end{enumerate}
\end{remark}

We conclude this section with the main result.

\begin{theorem}
 Assume that $s > 0$ and $f_t \in L^{\max\{s,p\}}([-A,A]^d)$ for some $p>1$. Let $X$ be an infinitely divisible random field and the control measure $\lambda$ of the infinitely divisible random measure be the Lebesgue measure. Let $\varepsilon > 0$. If $\tilde{f}_t^{(n)}$ is calculated using the algorithm mentioned above, then
 $$Err_s(X(t),\tilde{X}^{(n)}(t)) \leq \varepsilon, \quad \forall t \in [-T,T]^q.$$
\end{theorem}

\section{Simulation study}
\label{sec:simstud}

\begin{figure}[!ht]
\begin{center}
   \includegraphics[width=7cm]{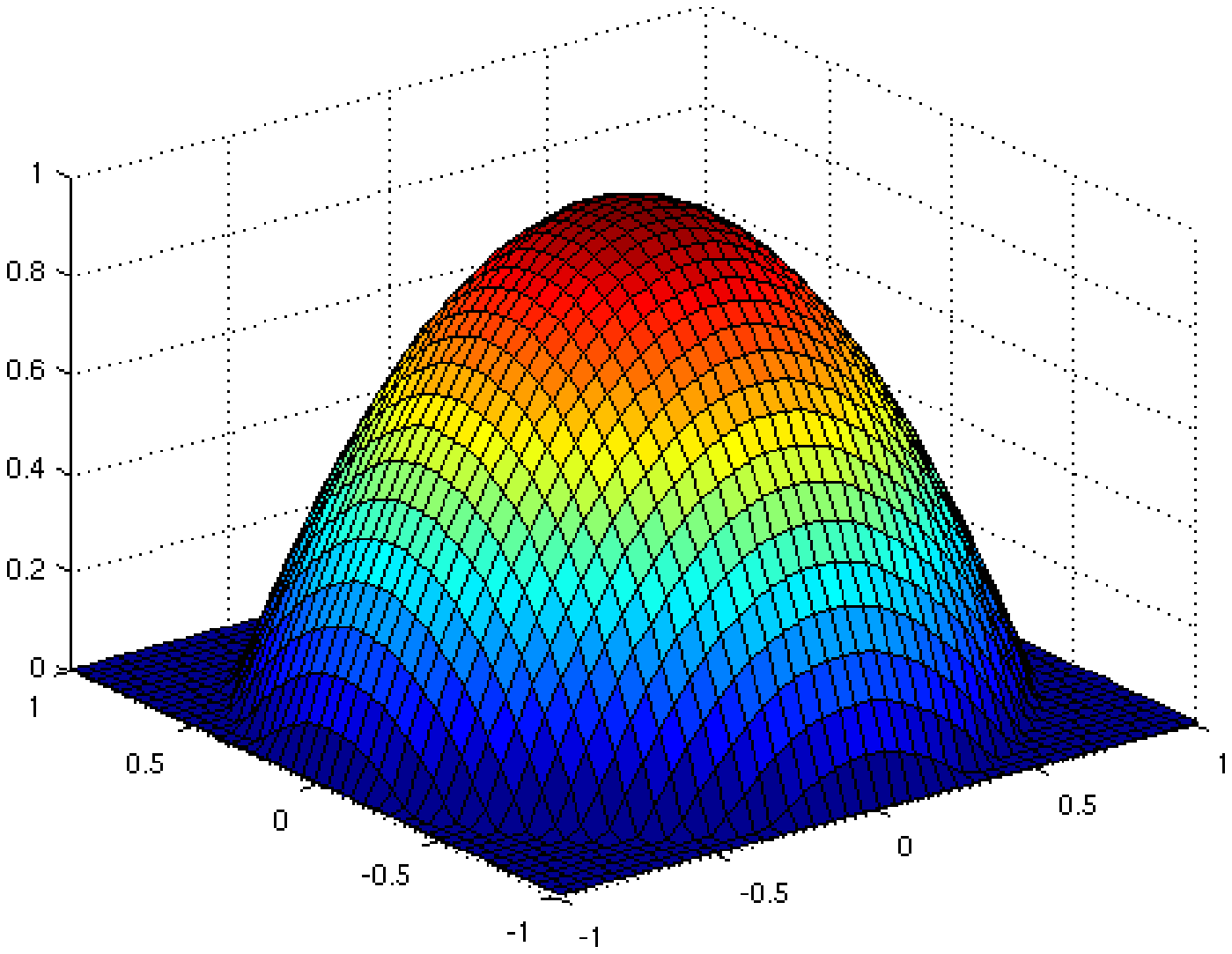}\includegraphics[width=7cm]{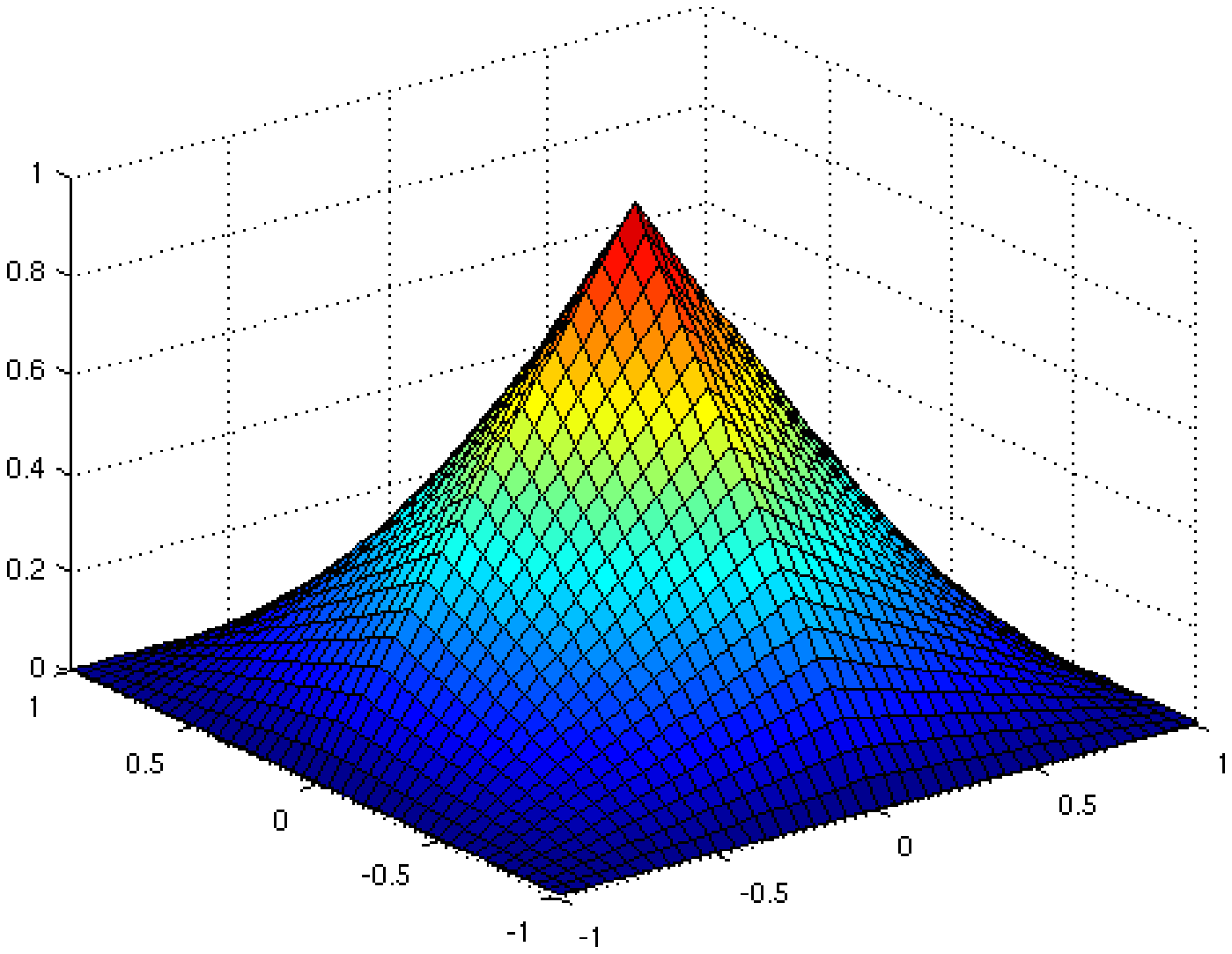}
   \end{center}
   \caption{The Epanechnikov-type kernel function (\ref{eq:epanechnikov}) (left) and the kernel function (\ref{eq:pyramid}) (right).}
   \label{fig:epanechnikov_pyramid}
 \end{figure}

For the simulation study, we used two different types of kernel functions for $\alpha$-stable random fields of dimension $d=2$. The first one is an Epanechnikov-type kernel function defined by
\begin{equation}
 f_{t}(x) = \begin{cases}
                     b \cdot (a^2-\Vert x-t \Vert_2^2), & \Vert x-t \Vert_2 \leq a, \\
		     0, & \text{otherwise},
                  \end{cases} \label{eq:epanechnikov}
\end{equation}
where $a>0$ and $b>0$, whereas for the second one, we take
\begin{equation} \label{eq:pyramid}
 \hspace{-0.4cm}\begin{split}
  f_{(t_1,t_2)}(x_1,x_2)& = b (a-|x_1 - t_1|) (a-|x_2-t_2|) \\
	& \hspace{0.3cm}\cdot\ind_{\{a-|x_1 - t_1| \geq 0,\text{\ } a-|x_2-t_2| \geq 0\}}(x_1,x_2)
 \end{split}
\end{equation}
where $a>0$ and $b>0$. Examples of both kernel functions are plotted in Figure \ref{fig:epanechnikov_pyramid}.

The main difference between these two types of kernel functions is that one can derive a simple formula for the integral of (\ref{eq:pyramid}) over squares, but not for the integral of (\ref{eq:epanechnikov}). This does not affect the step function approach since the kernel functions are only evaluated there at the points $\xi_k$, but it does affect the wavelet approach because the input vector consists of such integrals of (\ref{eq:epanechnikov}) and (\ref{eq:pyramid}) over squares. Therefore, we have to expect a loss in computational performance for kernel (\ref{eq:epanechnikov}) with the wavelet approach in this case.

Both functions (\ref{eq:epanechnikov}) and (\ref{eq:pyramid}) are H\"older-continuous with parameters $(C_1,\gamma_1) = (2ab,1)$ and $(C_2,\gamma_2) = (\sqrt{2}ab,1)$, respectively. We fixed $\alpha=1.5$, $\beta = 0$ and $[-T,T]^2=[-1,1]^2$ for both types of kernel functions.

\begin{figure}[ht]
   \begin{center}
   \includegraphics[height=5cm]{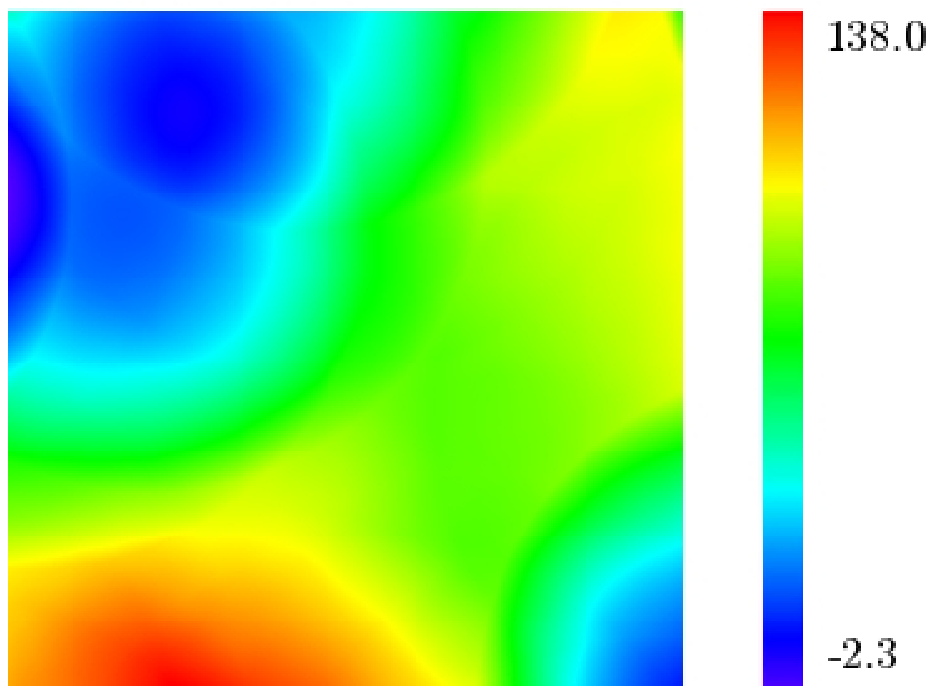}\hspace*{0.5cm}\includegraphics[height=5cm]{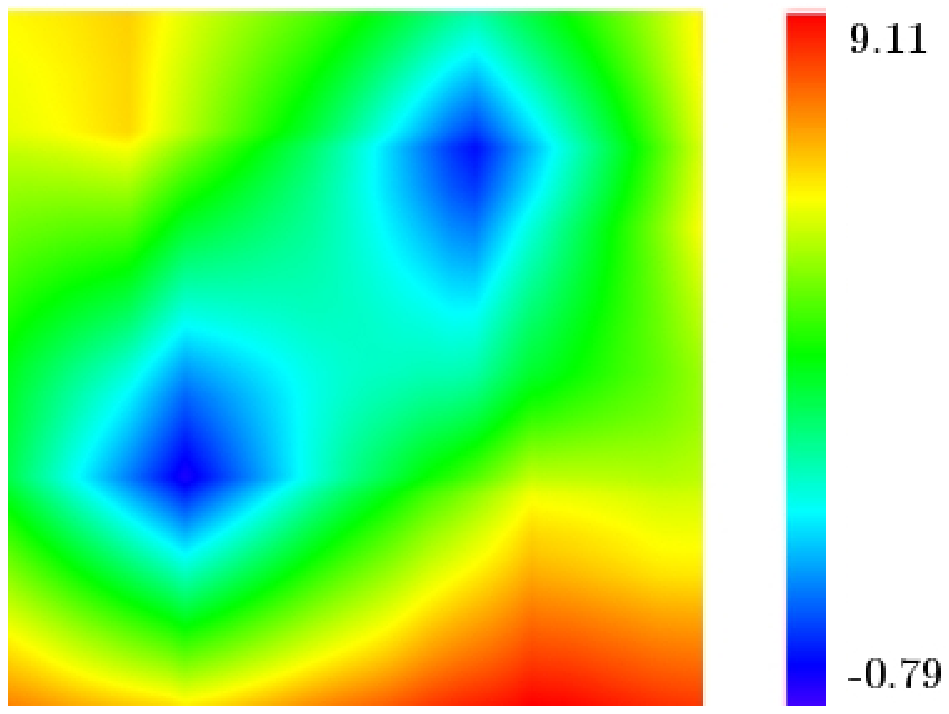}
   \end{center}
   \caption{Two realisations of stable random fields with kernel (\ref{eq:epanechnikov}) (left) and kernel (\ref{eq:pyramid}) (right).}
   \label{fig:epanechnikov_pyramid_realisation}
 \end{figure}

For the remaining parameters, we started with the following configuration: $b=1$, $a=1$ and $\varepsilon=1$. Furthermore, we divided $[-1,1]^2$ into an equidistant grid of $50 \times 50$ points and chose $l=0$ for the number of detail levels to be increased. Two realisations of the $1$-stable random field $X$ with kernels (\ref{eq:epanechnikov}) and (\ref{eq:pyramid}) are shown in Figure \ref{fig:epanechnikov_pyramid_realisation}.

First, we kept all parameters fixed and determined the computational time depending on the number of realisations. For the step function approach, each realisation needs the same computational time. For the wavelet approach, however, the wavelet coefficients only have to be calculated for the first realisation and can be stored afterwards. Therefore, any further realisation needs less computational time. Table \ref{tab:first_further} shows the results for both the Epanechnikov-type kernel (\ref{eq:epanechnikov}) and kernel (\ref{eq:pyramid}). By trial-and-error, we figured out that a combination of $\varepsilon=\varepsilon_1+\varepsilon_2 = 0.99 + 0.01$ performs quite good for the corresponding parameters in the wavelet algorithm in this case.

\begin{table}[htbp]
	\newcolumntype{d}[1]{D{.}{.}{#1}}
	\setlength{\tabcolsep}{10pt}
	\footnotesize
	\vspace{1.5mm}
	\caption{Computational time (in msec) for the first and further realisations.}
        \begin{tabular}{lcc}
	    \toprule
		& Kernel (\ref{eq:epanechnikov}) & Kernel (\ref{eq:pyramid}) \\
	    \midrule
		Step function & 28.5 & 18.6 \\
		approach & & \\
		\addlinespace
		Wavelet approach & 5337.0 & 1090.0 \\
		(first realisation) & & \\
		\addlinespace
		Wavelet approach & 13.5 & 13.2 \\
		(further realisations) & & \\
            \bottomrule
        \end{tabular}\\
	\label{tab:first_further}
\end{table}

Second, we focused on kernel (\ref{eq:pyramid}) and the computational time of any further realisation except the first one and varied subsequently one of the parameters $\alpha$, $m$ (the number of pixels per row) and $\varepsilon$ while all the other parameters were kept fixed. It turned out that the computational time decreased for the wavelet approach and increased for the step function approach when $\alpha$ was decreased. The computational time was equal for both approaches at about $\alpha=1.8$. For decreasing $\varepsilon$, the computational time for the step function approach increased much faster than the one for the wavelet approach. Varying $m$ affected the computational time of both approaches in a similar manner. The above results imply that neither of the two approaches outperforms the other one. For some combinations of the parameters, the step function approach was faster than the wavelet appraoch, for others it was slower.

Finally, we increased the parameter $l$ successively for a field with $10\times 10$ pixels while all other parameters were kept fixed and investigated the computational time for the wavelet approach for any further realisation except the first one. Table \ref{tab:l} shows the corresponding results.

\begin{table}[htbp]
	\newcolumntype{d}[1]{D{.}{.}{#1}}
	\footnotesize
	\vspace{1.5mm}	
	\caption{Computational time (in msec) for different values of $l$ (kernel (\ref{eq:pyramid})).}
        \begin{tabular}{cccccccc}
	    \toprule
		$l$ & 0 & 1 & 2 & 3 & 4 \\
	    \midrule
		Computational time & 25.5 & 45.5 & 246.4 & 1044.8 & 4212.0 \\
            \bottomrule
        \end{tabular}\\
	\label{tab:l}
\end{table}

One might have expected that the computational time tends to decrease if $l$ is increased since the wavelet series usually consists of less summands when keeping the same level of precision. At the same time, however, more stable random variable simulations have to be performed for the calculation of the integrals 
	$$ \int_{[-A,A]^d} \Psi_{j-2^{m+l},m+l}^e M(dx).$$
That is why for larger values of $l$, the computational time increases sharply.

For $\alpha = 2$, the random measure $M$ is a Gaussian random measure and the random field

\begin{equation}
 X(t) = \int_{\R^2} f_t(x) M(dx) \label{eq:stoch_int}
\end{equation}

becomes a Gaussian random field. Therefore, we can use the simulation methods to simulate such fields, too.

The kernel function

$$ f_t(x) = b \cdot e^{-\frac{\Vert x-t \Vert_2^2}{a}}, \quad t \in \R^2$$

corresponds to the isotropic and stationary covariance function

\begin{equation}
 C(h) := \Cov(X(0,0),X(h,0)) = 2 \int_{\R^2} b e^{-\frac{\Vert x \Vert_2^2}{a}} \cdot b e^{-\frac{\Vert x - (h,0)^T \Vert_2^2}{a}} dx = \pi a b^2 e^{-\frac{h^2}{2a}}, \label{eq:cov}
\end{equation}

cf.~\cite{ST94}, p.~128.

We now want to investigate how our simulation methods perform compared to the circulant embedding method (see~\cite{WC94}) to simulate Gaussian random fields with the given covariance function (\ref{eq:cov}) for $a=0.05$ and $b=1$ and $[-T,T]^2=[-0.5,0.5]^2$ for a grid of $100 \times 100$ points. We choose the circulant embedding method since it is exact in principle, and if exact simulation takes too much computational time, approximation techniques exist such that at least the one-dimensional marginal distributions are exact in principle.

For the comparison of the circulant embedding method with the step function approach and the wavelet approach, we simulate 1000 fields and estimate their mean, their variance and their covariance function for the distances 0, 0.01, 0.02, ..., 0.5. We then compare the values to the theoretical ones and require that the estimated values do not differ from the theoretical ones more than 0.01. If at least one value differs more than 0.01, the precision level is increased. The following table shows the computational time for each of the three methods.

\begin{table}[htbp]
	\newcolumntype{d}[1]{D{.}{.}{#1}}
	\footnotesize
	\vspace{1.5mm}	
	\caption{Computational time (in msec) for the circulant embedding method and the step function and the wavelet approach.}
        \begin{tabular}{cccc}
	    \toprule
		& Circulant embedding & Step function approach & Wavelet approach \\
	    \midrule
		Computational time & 48.43 & 9588.29 & 505.28 \\
            \bottomrule
        \end{tabular}\\
	\label{tab:circulantembedding}
\end{table}

As one can see, the wavelet approach and the step function approach take much more computational time than the circulant embedding method. This comes from the extensive calculations in the numerical integration of the stochastic integral (\ref{eq:stoch_int}). For Gaussian random fields, it is therefore advisable to use existing simulation methods such as the circulant embedding method that exploit the specific structure of these fields.

\section{Summary}
\label{sec:summary}

We presented two approaches to simulate $\alpha$-stable random fields that are based on approximating the kernel function by a step function and by a wavelet series. For both approaches, we derived estimates for the approximation error $Err(X(t), \tilde{X}^{(n)}(t))$.

In the simulation study we saw that for the first realisation of an $\alpha$-stable random field, the step function approach performs better than the wavelet approach due to the initial calculation of the wavelet coefficients. For any further realisation, however, the wavelet approach outperforms the step function approach.

Let us compare the rates of convergence of the step function approach and the wavelet approach more generally. If $f_t$ is H\"older-continuous, the error estimate for the step function approach is
\begin{equation}
 Err_s(X(t),\tilde{X}^{(n)}(t)) \leq C_1(d,C_t,\gamma_t,s,A) \cdot \left(\frac{1}{n}\right)^{\gamma_t} \label{eq:estimate_step}
\end{equation}
for a constant $C_1(d,C_t,\gamma_t,s,A)>0$. In the simulation study, we have seen that increasing the parameter $l$ in order to get closer to the best $n$-term approximation is not so advantageous. Therefore, we consider the rate of convergence for the cut wavelet series with error estimate
\begin{equation}
 Err_s(X(t),\tilde{X}^{(n)}(t)) \leq C_2(d,C_t,\gamma_t,s,A) \cdot \left(\frac{1}{2^{d/s+\gamma_t}}\right)^n \label{eq:estimate_wavelet}
\end{equation}
with a constant $C_2(d,C_t,\gamma_t,s,A)>0$. We note that we cannot compare the error estimates directly because for the step function approach, $n$ determines the number of cubes ($(2n)^d$) that form a partition of $[-A,A]^d$, while for the wavelet approach, $n$ is the detail level. Therefore, we express the error bounds in terms of the number of summands of the step function approximation (\ref{eq:rf_step_function}) of the random field, p.~\pageref{eq:rf_step_function}, and the wavelet approximation (\ref{eq:rf_wavelet}), p.~\pageref{eq:rf_wavelet}, respectively.
\begin{table}[htbp]
	\newcolumntype{d}[1]{D{.}{.}{#1}}
	\footnotesize
	\vspace{1.5mm}
	\caption{Number of summands of (\ref{eq:rf_step_function}) and (\ref{eq:rf_wavelet}) and error bounds for (\ref{eq:estimate_step}) and (\ref{eq:estimate_wavelet}) in terms of the number of summands.}
        \begin{tabular}{ccc}
	    \toprule
		& Step function approach & Wavelet approach \\
	    \midrule
	    \addlinespace
		Number of summands & $u=(2n)^d$ & $u=1+2^{d-1}d(2^{n+1}-1)$ \\
		Error bounds & $O\left(\left(\frac{1}{u}\right)^{\gamma_t/d}\right)$ & $O\left(\left(\frac{1}{2^{(d/s+\gamma_t)/\ln(2)}}\right)^{\ln\left(\frac{2(u-1)}{2^d d}+1\right)}\right)$\\
	    \addlinespace
            \bottomrule
        \end{tabular}\\
	\label{tab:rates}
\end{table}
In Table \ref{tab:rates}, we see that the rates of convergence in terms of the number of summands do not distinguish substantially between the step function approach and the wavelet approach.

Let us consider the number of random variables that need to be simulated for a single realisation of a random field. For the step function approach, we need to simulate $(2n)^d$ random variables. For the wavelet approach, the number of random variables to simulate is equal to the number of cubes that form a partition of $[-A,A]^d$ in the finest detail level $n$: $2^{d(n+1)}$. In the examples of the simulation study, much less random variables had to be simulated for the wavelet approach than for the step function approach which was a reason for the good performance of the wavelet approach.

We want to make two remarks about the wavelet approach. First, we have seen that one drawback is that the computation of the input vector for the fast wavelet transfrom may take quite a long time if no formula for the integrals
$$ \int_{C} f_t(x) dx $$
is known, where $C$ is a cube in $\R^d$. In general, for an arbitrary wavelet basis $\{\Psi_i^*\}_{i \in I}$, we would have to calculate
$$ \int_C f_t(x) \Psi_i(x) dx$$
for some $i \in I$ which, in many cases, also requires numerical integration.

Interpolatory wavelet bases can remedy this disadvantage since for this kind of wavelets bases, the wavelet coefficients basically reduce to evaluating the kernel function at a certain point. However, the interpolatory wavelets themselves are no step functions any more such that the simulation of the integrals
$$ \int_{[-A,A]^d} \Psi_i(x) \Lambda(dx),$$
where $\Psi_i$ is an interpolatory wavelet function, is much more complicated than for the Haar basis.

Second, one could use adaptive wavelet methods in order to calculate the wavelet coefficients. This might decrease the computational time for the first random field realisation. However, in the simulation study we have seen that increasing the parameter $l$ has little advantage over the cut wavelet series ($l=0$) since the negative effect of the increasing detail level and thus the need of more stable random variable simulations dominates the positive one of less summands in the wavelet series.

\textbf{Acknowledgements.} The authors wish to thank Prof. Urban for his assistance in wavelet-related questions. They also want to thank Generali Versicherung AG, Vienna, Austria, for the kind support of this research.

\end{document}